\documentclass[aip,jcp, reprint,groupedaddress]{revtex4-2}%
\usepackage{graphicx}
\usepackage{tabularx}
\usepackage{dcolumn}
\usepackage{longtable}
\usepackage{tensor}
\usepackage{color}
\usepackage{placeins}
\usepackage{bm}
\usepackage{amsmath}
\usepackage{amsfonts}
\usepackage{amssymb}
\usepackage{float}
\usepackage[urlcolor=white]{hyperref}%
\providecommand{\U}[1]{\protect\rule{.1in}{.1in}}
\graphicspath{{C:/Users/moghadda/Desktop/VGALCA/Figures/}{C:/Users/Jiri/Dropbox/Papers/Chemistry_papers/2023/VGA-LCA/Figures/}}
\begin{document}
\title{ High-order geometric integrators for the local cubic variational Gaussian
wavepacket dynamics}
\author{Roya Moghaddasi Fereidani}
\email{roya.moghaddasifereidani@epfl.ch}
\author{Ji\v{r}\'i J. L. Van\'i\v{c}ek}
\email{jiri.vanicek@epfl.ch}
\affiliation{Laboratory of Theoretical Physical Chemistry, Institut des Sciences et
Ing\'enierie Chimiques, Ecole Polytechnique F\'ed\'erale de Lausanne (EPFL),
CH-1015, Lausanne, Switzerland}
\date{\today}

\begin{abstract}
Gaussian wavepacket dynamics has proven to be a useful semiclassical
approximation for quantum simulations of high-dimensional systems with low
anharmonicity. Compared to Heller's original local harmonic method, the
variational Gaussian wavepacket dynamics is more accurate, but much more
difficult to apply in practice because it requires evaluating the expectation
values of the potential energy, gradient, and Hessian. If the variational
approach is applied to the local cubic approximation of the potential, these
expectation values can be evaluated analytically, but still require the costly
third derivative of the potential. To reduce the cost of the resulting local
cubic variational Gaussian wavepacket dynamics, we describe efficient
high-order geometric integrators, which are symplectic, time-reversible, and
norm-conserving. For small time steps, they also conserve the effective
energy. We demonstrate the efficiency and geometric properties of these
integrators numerically on a multi-dimensional, nonseparable coupled Morse
potential.

\end{abstract}
\maketitle

\section{\label{sec:introduction}Introduction}

Gaussian wavepackets provide an explicit link between the classical and
quantum dynamics.~\cite{book_Tannor:2007} Because they are localized in space,
Gaussian wavepackets only require information about the potential in a
specific region of configuration space, and thus can substantially reduce the
cost of molecular quantum dynamics simulations. Inspired by the pioneering
work of Heller,~\cite{Heller:1975, Heller:1976a, Heller:1981} a number of
multi-trajectory Gaussian-based methods,\textbf{ }such as the Herman-Kluk
initial value
representation,~\cite{Herman_Kluk:1984,Walton_Manolopoulos:1995,Miller:2001}
multiple spawning,~\cite{Martinez_Levine:1996a,
Ben-Nun_Martinez:2000,Curchod_Martinez:2018_v2} variational
multiconfigurational Gaussians,\cite{Worth_Burghardt:2003} coupled coherent
states,\cite{Shalashilin_Child:2000} multiconfigurational Ehrenfest
method,\cite{Shalashilin:2009} or Gaussian dephasing
representation,~\cite{Sulc_Vanicek:2013} were proposed to approximate the
solution of the time-dependent Schr\"{o}dinger equation (TDSE). Because the
local nature of Gaussian wavepackets allows the potential to be calculated
\textquotedblleft on the fly,\textquotedblright\ many of these methods were
combined with an on-the-fly ab initio evaluation of the electronic
structure.~\cite{Tatchen_Pollak:2009,Ceotto_Atahan:2009a,Wong_Roy:2011,Botti_Conte:2021,Levine_Coe:2008,Levine_Martinez:2009,Bramley_Shalashilin:2019,Worth_Burghardt:2004,Saita_Shalashilin:2012}%

The multi-trajectory methods typically require a large number of trajectories
to converge, and therefore, single-trajectory Gaussian-based methods have
recently re-emerged as a practical alternative. Although they cannot capture
wavepacket splitting, the single-trajectory methods, such as
divide-and-conquer semiclassical
dynamics,\cite{Ceotto_Conte:2017,Gabas_Ceotto:2019} ab initio thawed Gaussian
approximation,~\cite{Wehrle_Vanicek:2014_v2,Wehrle_Vanicek:2015_v2,Begusic_Vanicek:2020a,Kletnieks_Vanicek:2023}
and its more efficient single Hessian
version,~\cite{Begusic_Vanicek:2019,Prlj_Vanicek:2020,Begusic_Vanicek:2022,Golubev_Vanicek:2020}
have proven sufficiently accurate in weakly anharmonic systems of rather large
dimension. Moreover, thanks to their efficiency, these methods can afford more
accurate electronic structure theory.

Among the single-trajectory Gaussian-based methods, the most accurate one is
the variational Gaussian wavepacket dynamics (variational
GWD),~\cite{Heller:1976,
Coalson_Karplus:1990,book_Lubich:2008,Lasser_Lubich:2020,Vanicek:2023,Fereidani_Vanicek:2023}
obtained by applying the time-dependent variational
principle~\cite{Dirac:1930,book_Frenkel:1934,book_Lubich:2008} to a thawed
Gaussian ansatz (\textquotedblleft thawed\textquotedblright\ means that the
Gaussian has a flexible width, which we will always assume). Not only is the
variational GWD symplectic and time-reversible, but it also conserves the norm
and energy of the wavepacket.\cite{Faou_Lubich:2006,Fereidani_Vanicek:2023}
The method is, however, hard to apply to \textquotedblleft
real-world\textquotedblright\ problems because it requires expectation values
of the potential and its derivatives, which typically cannot be evaluated
analytically and are also very difficult to compute numerically. One way to
surmount this obstacle is to approximate the potential $V$ by its Taylor
expansion $V_{\text{appr }}$around the center of the Gaussian since the
expectation values of $V_{\text{appr}}$ and its derivatives are readily
evaluated analytically. Replacing the expectation values of the original
potential, gradient, and Hessian with the corresponding values of
$V_{\text{appr}}$ in the equations of motion of the variational GWD yields a
method, which is easily applied to any potential $V$, including that obtained
from an on-the-fly ab initio evaluation of the electronic structure. Using the
local harmonic approximation of the potential ($V_{\text{appr}}=V_{\text{LHA}%
}$) reduces the variational method to Heller's celebrated thawed Gaussian
approximation,~\cite{Heller:1975,Lee_Heller:1982,book_Heller:2018} which is,
however, neither symplectic nor energy-conserving and cannot capture
tunneling.~\cite{Fereidani_Vanicek:2023}

In contrast, the local cubic variational GWD, obtained by applying the local
cubic approximation to the potential ($V_{\text{appr}}=V_{\text{LCA}}$), is
symplectic,~\cite{Ohsawa:2015a,Ohsawa:2015b} conserves the effective
energy,~\cite{Pattanayak_Schieve:1994a,Vanicek:2023} and may capture tunneling
because the guiding trajectory is nonclassical.~\cite{Ohsawa_Leok:2013} This
remarkable method was originally studied by Pattanayak and
Schieve~\cite{Pattanayak_Schieve:1994a, Pattanayak_Schieve:1994b} under the
name \textquotedblleft extended semiclassical wavepacket
dynamics\textquotedblright\ and later, in depth, by Ohsawa and
co-workers,~\cite{Ohsawa_Leok:2013,Ohsawa:2015a,Ohsawa:2015b} who called it
\textquotedblleft symplectic semiclassical wavepacket
dynamics.\textquotedblright\ The local cubic variational GWD, which is also
norm-conserving and time-reversible, preserves all the geometric properties of
the fully variational GWD, except that it conserves the effective rather than
the exact energy. The improvements over Heller's local harmonic GWD, however,
come at a high cost of evaluating the symmetric rank-three tensor of third
derivatives of $V$ along the trajectory.

To make this appealing method more practical, in this paper, our goal is
reducing the computational cost of the local cubic variational GWD without
sacrificing its geometric properties. In the case of local harmonic GWD, it
was shown that using a single, i.e., constant Hessian along the whole
trajectory speeds up the method and, in addition, results in the conservation
of both symplectic structure and effective
energy.~\cite{Begusic_Vanicek:2019,Vanicek:2023} Applying a similar trick to
the local cubic variational GWD, namely, using a single third derivative,
however, breaks the symplecticity and effective energy conservation.
Therefore, to improve the efficiency of the local cubic method, we implement
and analyze high-order geometric integrators, which are analogous to the
integrators we recently implemented for the fully variational
GWD~\cite{Fereidani_Vanicek:2023} and which are a special case of general
integrators for GWD described in Ref.~\onlinecite{Vanicek:2023}.

The rest of this paper is organized as follows: In Sec.~\ref{sec:theory}, we
review the variational GWD and derive equations of motion for the GWD when the
potential is approximated by either the local harmonic or local cubic
expansion. After discussing the remarkable geometric properties of the local
cubic variational GWD in Sec.~\ref{sec:geom_prop}, in
Sec.~\ref{sec:geometric_integrators} we describe efficient high-order
geometric integrators that preserve most of these properties exactly.
Section~\ref{sec:numerical_example} provides numerical examples that
demonstrate the improved accuracy and geometric properties of the local cubic
variational GWD compared to the original local harmonic method and that
confirm the fast convergence, increased efficiency, and geometric properties
of the high-order integrators in various many-dimensional anharmonic systems.
Section~\ref{sec:conclusion} concludes the paper.


\section{Local cubic variational Gaussian wavepacket dynamics
\label{sec:theory}}

\subsection{Gaussian wavepacket dynamics}

The GWD approximates the solution of the TDSE
\begin{equation}
i\hbar\,d|\psi_{t}\rangle/dt=\hat{H}|\psi_{t}\rangle\label{eq:TDSE}%
\end{equation}
with a complex Gaussian ansatz,
\begin{equation}
\psi_{t}(q)=\text{exp}\bigg[\frac{i}{\hbar}\bigg(\frac{1}{2}\,x^{T}\cdot
A_{t}\cdot x+p_{t}^{T}\cdot x+\gamma_{t}\bigg)\bigg], \label{eq:GWP}%
\end{equation}
where $x:=q-q_{t}$ is the difference position vector, $q_{t}$ and $p_{t}$ are
the $D$-dimensional real vectors denoting the position and momentum of the
Gaussian's center, $A_{t}$ is a $D\times D$ complex symmetric width matrix
with a positive-definite imaginary part, and $\gamma_{t}$ is a complex number
whose real part represents a phase and whose imaginary part ensures
normalization of $\psi_{t}$. In Eq.~(\ref{eq:TDSE}), $\hat{H}\equiv H(\hat
{q},\hat{p})$ is a Hamiltonian operator separable into a sum
\begin{equation}
\hat{H}=T(\hat{p})+V(\hat{q}) \label{eq:H}%
\end{equation}
of the kinetic energy operator $T(\hat{p})=\hat{p}^{T}\cdot m^{-1}\cdot\hat
{p}/2$, depending only on momentum $p$, and the potential energy operator
$V(\hat{q})$, depending only on position $q$. $m$ is the real-symmetric mass
matrix. It can be shown~\cite{Vanicek:2023} that the GWD for
Eq.~(\ref{eq:TDSE}) can be viewed as the exact solution of the nonlinear
TDSE~\cite{Roulet_Vanicek:2021a,Roulet_Vanicek:2021b}
\begin{equation}
i\hbar\,d|\psi_{t}\rangle/dt=\hat{H}_{\text{eff}}(\psi_{t})|\psi_{t}%
\rangle\label{eq:NLTDSE}%
\end{equation}
with an effective, state-dependent Hamiltonian operator
\begin{equation}
\hat{H}_{\text{eff}}(\psi_{t})=T(\hat{p})+V_{\text{eff}}(\hat{q};\psi_{t}),
\label{eq:Heff}%
\end{equation}
where $V_{\text{eff}}(\hat{q};\psi_{t})$ is an effective quadratic
potential~\cite{Vanicek:2023,Fereidani_Vanicek:2023}
\begin{equation}
V_{\text{eff}}(q;\psi_{t})=V_{0}+V_{1}^{T}\cdot x+x^{T}\cdot V_{2}\cdot x/2.
\label{eq:Veff}%
\end{equation}
The coefficients $V_{0}$, $V_{1}$, and $V_{2}$ may depend on the state
$\psi_{t}$ and are different for each of the GWD methods. In a general GWD,
the parameters of the Gaussian wavepacket~(\ref{eq:GWP}) satisfy the system
\begin{align}
\dot{q}_{t}  &  =m^{-1}\cdot p_{t},\label{eq:qEOM}\\
\dot{p}_{t}  &  =-V_{1},\label{eq:pEOM}\\
\dot{A}_{t}  &  =-A_{t}\cdot m^{-1}\cdot A_{t}-V_{2},\label{eq:AEOM}\\
\dot{\gamma}_{t}  &  =T(p_{t})-V_{0}+(i\hbar/2)\,\text{Tr}(m^{-1}\cdot A_{t})
\label{eq:gammaEOM}%
\end{align}
of first-order ordinary differential equations.~\cite{Vanicek:2023}


\subsection{Fully variational GWD}

The variational GWD employs the time-dependent variational
principle~\cite{Dirac:1930,book_Frenkel:1934,book_Lubich:2008} to optimally
approximate the solution of the TDSE~(\ref{eq:TDSE}) with the Gaussian
ansatz~(\ref{eq:GWP}). In the variational GWD, the effective
potential~(\ref{eq:Veff}) has the
coefficients~\cite{Vanicek:2023,Fereidani_Vanicek:2023}
\begin{equation}
V_{0}=\langle\hat{V}\rangle-\text{Tr}\big({\langle\hat{V}^{\prime\prime
}\rangle}\cdot\Sigma_{t}\big)/2,\quad V_{1}=\langle\hat{V}^{\prime}%
\rangle,\quad V_{2}=\langle\hat{V}^{\prime\prime}\rangle.
\label{eq:V0_V1_V2_VGA}%
\end{equation}
Above and throughout this paper, we use a shorthand notation $\langle\hat
{O}\rangle:=\langle\psi_{t}|\hat{O}|\psi_{t}\rangle$ for the expectation value
of the operator $\hat{O}$ in the Gaussian wavepacket $\psi_{t}$. In
Eq.~(\ref{eq:V0_V1_V2_VGA}), $\hat{V}^{\prime}:=V^{\prime}(q)|_{q=\hat{q}}$
and $\hat{V}^{\prime\prime}:=V^{\prime\prime}(q)|_{q=\hat{q}}$ are the
gradient and Hessian of the potential energy operator, and
\begin{equation}
\Sigma_{t}:=\langle\hat{x}\otimes\hat{x}^{T}\rangle=(\hbar/2)\,(\text{Im}%
A_{t})^{-1} \label{eq:cov_q}%
\end{equation}
is the position covariance.~\cite{Vanicek_Begusic:2021} Therefore, the
equations of motion for the variational method include expectation values of
the potential and its derivatives, which generally cannot be evaluated
analytically. Since the expectation value of a polynomial operator in a
Gaussian can be obtained exactly, in the following sections, we approximate
the potential by its second- and third-order Taylor expansions about the
Gaussian's center and obtain the modified equations of motion for the
variational GWD.


\subsection{Local harmonic (variational) GWD = Heller's GWD}

Replacing the exact potential $V(q)$ in Eq.~(\ref{eq:V0_V1_V2_VGA}) with its
local harmonic approximation (LHA)
\begin{equation}
V_{\text{LHA}}(q,\psi_{t})=V(q_{t})+V^{\prime}(q_{t})^{T}\cdot x+x^{T}\cdot
V^{\prime\prime}(q_{t})\cdot x/2 \label{eq:V_LHA}%
\end{equation}
gives\cite{Vanicek:2023}
\begin{equation}
V_{0}=V(q_{t}),\quad V_{1}=V^{\prime}(q_{t}),\quad V_{2}=V^{\prime\prime
}(q_{t}). \label{eq:V0_V1_V2_TGA}%
\end{equation}
Equations~(\ref{eq:qEOM})-(\ref{eq:gammaEOM}) with
coefficients~(\ref{eq:V0_V1_V2_TGA}) are exactly the equations of motion of
Heller's original thawed Gaussian approximation.\cite{Heller:1975,
Wehrle_Vanicek:2014_v2} This local harmonic GWD is very efficient and easy to
implement, but preserves neither the effective nor the exact
energy.~\cite{Vanicek:2023}
Moreover, since the Gaussian's center $(q_{t},p_{t})$ evolves according to
classical Hamilton's equations of motion [Eqs.~(\ref{eq:qEOM}) and
(\ref{eq:pEOM}) with coefficients~(\ref{eq:V0_V1_V2_TGA})], the local harmonic
GWD cannot capture quantum tunneling.


\subsection{Local cubic variational GWD}

To improve upon the local harmonic GWD, one needs to replace the potential
with a more accurate approximation. Let us, therefore, apply the variational
GWD to the local cubic approximation (LCA)
\begin{align}
V_{\text{LCA}}(q,\psi_{t}) &  =V(q_{t})+V^{\prime}(q_{t})^{T}\cdot
x+x^{T}\cdot V^{\prime\prime}(q_{t})\cdot x/2\nonumber\\
&  \quad+V^{\prime\prime\prime}(q_{t})_{ijk}\,x_{i}x_{j}x_{k}%
/3!\label{eq:V_LCA}%
\end{align}
of the potential $V$, where $V^{\prime\prime\prime}$ is the rank-3 tensor of
the third derivatives of the potential $V$. The gradient and Hessian of the
$V_{\text{LCA}}$ are
\begin{align}
V_{\text{LCA}}^{\prime}(q,\psi_{t})_{i} &  =V^{\prime}(q_{t})_{i}%
+V^{\prime\prime}(q_{t})_{ij}\,x_{j}\nonumber\\
&  \quad+V^{\prime\prime\prime}(q_{t})_{ijk}\,x_{j}x_{k}%
/2,\label{eq:grad_V_LCA}\\
V_{\text{LCA}}^{\prime\prime}(q,\psi_{t})_{ij} &  =V^{\prime\prime}%
(q_{t})_{ij}+V^{\prime\prime\prime}(q_{t})_{ijk}\,x_{k}.\label{eq:Hess_V_LCA}%
\end{align}
By substituting expectation values
\begin{align}
\langle\hat{V}_{\text{LCA}}\rangle &  =V(q_{t})+\text{Tr}\big[V^{\prime\prime
}(q_{t})\cdot\Sigma_{t}\big]/2,\nonumber\\
\langle\hat{V}_{\text{LCA}}^{\prime}\rangle_{i} &  =V^{\prime}(q_{t}%
)_{i}+V^{\prime\prime\prime}(q_{t})_{ijk}\,\Sigma_{t,jk}%
/2,\label{eq:exp_value_Hess_V_LCA}\\
\langle\hat{V}_{\text{LCA}}^{\prime\prime}\rangle &  =V^{\prime\prime}%
(q_{t})\nonumber
\end{align}
into Eq.~(\ref{eq:V0_V1_V2_VGA}), we get the effective potential coefficients
\begin{align}
V_{0} &  =V(q_{t}),\nonumber\\
V_{1,i} &  =V^{\prime}(q_{t})_{i}+V^{\prime\prime\prime}(q_{t})_{ijk}%
\,\Sigma_{t,jk}/2,\label{eq:V0_V1_V2_LCA}\\
V_{2} &  =V^{\prime\prime}(q_{t}).\nonumber
\end{align}
The equations of motion of the resulting \textquotedblleft local cubic
variational GWD\textquotedblright,
\begin{align}
\dot{q}_{t} &  =m^{-1}\cdot p_{t},\label{eq:qEOM_LCA}\\
\big(\dot{p}_{t}\big)_{i} &  =-V^{\prime}(q_{t})_{i}-V^{\prime\prime\prime
}(q_{t})_{ijk}\,\Sigma_{t,jk}/2,\label{eq:pEOM_LCA}\\
\dot{A}_{t} &  =-A_{t}\cdot m^{-1}\cdot A_{t}-V^{\prime\prime}(q_{t}%
),\label{eq:AEOM_LCA}\\
\dot{\gamma}_{t} &  =T(p_{t})-V(q_{t})+(i\hbar/2)\,\text{Tr}(m^{-1}\cdot
A_{t}),\label{eq:gammaEOM_LCA}%
\end{align}
are identical to those of the local harmonic GWD except for a non-classical
force term $-V^{\prime\prime\prime}(q_{t})_{ijk}\,\Sigma_{t,jk}/2$ in the
equation for $\dot{p}_{t}$. The local cubic variational GWD is equivalent to
Ohsawa and Leok's \textquotedblleft symplectic semiclassical wavepacket
dynamics\textquotedblright~\cite{Ohsawa_Leok:2013,Ohsawa:2015a,Ohsawa:2015b}
and Pattanayak and Schieve's \textquotedblleft extended semiclassical
dynamics\textquotedblright.~\cite{Pattanayak_Schieve:1994b} Remarkably, in
contrast to the local harmonic GWD, the local cubic method is symplectic,
conserves the effective energy,~\cite{Vanicek:2023} and may capture tunneling
because the Gaussian follows a nonclassical
trajectory.~\cite{Ohsawa_Leok:2013}

Let us mention that the utility of the local cubic expansion of the PES is not
restricted to single-trajectory methods---this expansion was also used for
evaluating the matrix elements of the potential needed in the calculation of
the absorption spectrum of formaldehyde with the on-the-fly variational
multiconfigurational Gaussian method; there, too, the local cubic expansion
was found to yield more accurate results than the local harmonic
approximation.~\cite{Bonfanti_Pollak:2018}

\section{Geometric properties of the local cubic variational GWD
\label{sec:geom_prop}}

In the local cubic variational GWD, the Gaussian wavepacket evolves according
to a nonlinear TDSE~(\ref{eq:NLTDSE}) in a state-dependent effective quadratic
potential (\ref{eq:Veff}) with coefficients~(\ref{eq:V0_V1_V2_LCA}). As any
other exact solution of the more general nonlinear TDSE~(\ref{eq:NLTDSE}), the
local cubic variational GWD is time-reversible and norm-conserving, but does
not conserve the inner product and
distance.~\cite{Roulet_Vanicek:2021a,Vanicek:2023} Unlike some other GWD
methods, the local cubic variational method is
symplectic~\cite{Ohsawa_Leok:2013} and, although it is not energy-conserving,
it does conserve the \emph{effective} energy exactly.~\cite{Vanicek:2023} Let
us discuss the last three properties in more detail.


\subsection{Energy}

The energy of a Gaussian wavepacket can be computed as the expectation value
\begin{equation}
E:=\langle\hat{H}\rangle=\langle\hat{T}\rangle+\langle\hat{V}\rangle
\label{eq:E_exact}%
\end{equation}
of the separable Hamiltonian~(\ref{eq:H}). In Eq.~(\ref{eq:E_exact}),
\begin{equation}
\langle\hat{T}\rangle=T(p_{t})+\text{Tr}\big[m^{-1}\cdot\text{Cov}(\hat
{p})\big]/2 \label{eq:E_kin}%
\end{equation}
is the kinetic energy,~\cite{Vanicek:2023} $\langle\hat{V}\rangle$ is the
potential energy, and
\begin{align}
\text{Cov}(\hat{p})  &  :=\langle(\hat{p}-p_{t})\otimes(\hat{p}-p_{t}%
)^{T}\rangle\nonumber\\
&  \,\,=(\hbar/2)\,A_{t}\cdot(\text{Im}\,A_{t})^{-1}\cdot A_{t}^{\dag}
\label{eq:Covp}%
\end{align}
is the momentum covariance.~\cite{Vanicek:2023} In general GWD, the energy
depends on time as [see Eq. (21) in Ref.~\onlinecite{Vanicek:2023}]:
\begin{equation}
\dot{E}=\text{Re}\,\big \langle\hat{p}^{T}\cdot m^{-1}\cdot(\hat{V}^{\prime
}-\hat{V}_{\text{eff}}^{\prime})\big \rangle.
\end{equation}
For the local cubic method with coefficients~(\ref{eq:V0_V1_V2_LCA}), we have
\begin{align}
V_{\text{eff}}^{\prime}(q,\psi_{t})_{i}  &  =V^{\prime}(q_{t})_{i}%
+V^{\prime\prime\prime}(q_{t})_{ijk}\,\Sigma_{t,jk}/2+V^{\prime\prime}%
(q_{t})_{ij}\,x_{j}\nonumber\\
&  \neq V^{\prime}(q)_{i}.
\end{align}
Therefore, the local cubic variational GWD does not conserve
energy,\cite{Vanicek:2023} even though the fully variational GWD is
energy-conserving.~\cite{Fereidani_Vanicek:2023, Lasser_Lubich:2020}

\subsection{Effective energy}

The effective energy of a Gaussian wavepacket driven by the effective
Hamiltonian~(\ref{eq:Heff}) can be evaluated as
\begin{equation}
E_{\text{eff}}:=\langle\hat{H}_{\text{eff}}\rangle=\langle\hat{T}%
\rangle+\langle\hat{V}_{\text{eff}}\rangle, \label{eq:E_eff}%
\end{equation}
where $\langle\hat{T}\rangle$ is given in Eq.~(\ref{eq:E_kin}) and
$\langle\hat{V}_{\text{eff}}\rangle$ for the local cubic variational GWD is
obtained from Eqs.~(\ref{eq:Veff}) and~(\ref{eq:V0_V1_V2_LCA}),
\begin{equation}
\langle\hat{V}_{\text{eff}}\rangle=V(q_{t})+\text{Tr}\big[V^{\prime\prime
}(q_{t})\cdot\Sigma_{t}\big]/2. \label{eq:Veff_exp_value_LCA}%
\end{equation}
The time dependence of the effective energy~(\ref{eq:E_eff}) is [see Eq. (22)
in Ref.~\onlinecite{Vanicek:2023}]:
\begin{equation}
\dot{E}_{\text{eff}}=\dot{V}_{0}-V_{1}^{T}\cdot\dot{q}_{t}+\text{Tr}%
\big(\dot{V}_{2}\cdot\Sigma_{t}\big)/2. \label{eq:E_eff_conservation_3}%
\end{equation}
Inserting \textquotedblleft local cubic\textquotedblright%
\ coefficients~(\ref{eq:V0_V1_V2_LCA}) into Eq.~(\ref{eq:E_eff_conservation_3}%
) gives
\begin{align}
\dot{E}_{\text{eff}}  &  =V^{\prime}(q_{t})^{T}\cdot\dot{q}_{t}-V^{\prime
}(q_{t})^{T}\cdot\dot{q}_{t}-V^{\prime\prime\prime}(q_{t})_{ijk}%
\,\Sigma_{t,jk}\dot{q}_{t,i}/2\nonumber\\
&  \quad+V^{\prime\prime\prime}(q_{t})_{ijk}\,\Sigma_{t,jk}\dot{q}_{t,i}/2=0,
\label{eq:E_eff_conservation_4}%
\end{align}
and thus, unlike the local harmonic GWD, the local cubic variational method
does conserve the effective energy.

\subsection{Symplecticity}

Faou and Lubich demonstrated the non-canonical symplectic structure of the
variational GWD.~\cite{Faou_Lubich:2006} Ohsawa and Leok showed that the local
cubic variational GWD is a non-canonical Hamiltonian system whose equations of
motion are precisely Eqs.~(\ref{eq:qEOM_LCA})-(\ref{eq:gammaEOM_LCA}%
).~\cite{Ohsawa_Leok:2013} Omitting the phase factor $\exp(i\,\text{Re}%
\,\gamma_{t}/\hbar)$ in (\ref{eq:GWP}) corresponds to an equivalent
representation of the Gaussian wavepacket (without a phase) in a projective
Hilbert space with a phase symmetry.~\cite{Ohsawa:2015b, Ohsawa:2015a}
Furthermore, the imaginary part of the parameter $\gamma_{t}$ ensures
normalization of the Gaussian wavepacket and can be expressed in terms of the
imaginary part of the width $A_{t}$ [Eq.~(\ref{eq:Norm})]. Without the
parameter $\gamma$, the symplectic manifold parametrized by the set
$(q,p,\text{Re}\,A,\text{Im}\,A)$ has the reduced symplectic
form~\cite{Ohsawa:2015b, Ohsawa:2015a}
\begin{equation}
\omega(z)=dq_{j}\wedge dp_{j}+(\hbar/4)\,d\,\big[(\text{Im}\,A)^{-1}%
\big]_{kl}\wedge d\,\text{Re}\,A_{kl},\label{eq:sympl_Agamma}%
\end{equation}
where the first term is the standard classical symplectic structure and the
second term could be called the \textquotedblleft quantum\textquotedblright%
\ component of the symplectic structure.

\section{Geometric integrators for the local cubic variational GWD
\label{sec:geometric_integrators}}

Geometric integrators for the nonlinear TDSE~(\ref{eq:NLTDSE}) with a general
effective potential~(\ref{eq:Veff}) were presented in
Ref.~\onlinecite{Vanicek:2023}. The only and rather weak assumption necessary
for obtaining explicit integrators was that the coefficients $V_{0}$, $V_{1}$,
and $V_{2}$ of the effective potential~(\ref{eq:Veff}) depended on the
wavepacket $\psi_{t}$ only through $q_{t}$ and $\text{Im}\,A_{t}$ and were
independent of $p_{t}$, $\text{Re}\,A_{t}$, and $\gamma_{t}$. In
Ref.~\onlinecite{Fereidani_Vanicek:2023}, we investigated these geometric
integrators for the variational GWD and showed that the high-order integrators
could be more efficient and more accurate at the same time. Here, we describe,
implement, and analyze such geometric integrators for the local cubic
variational GWD.

\subsection{Second-order geometric integrator}

The simplest, first-order algorithm is based on the decomposition of the
separable, effective Hamiltonian~(\ref{eq:Heff}) into the kinetic and
potential terms. This makes it possible to solve the equations of
motion~(\ref{eq:qEOM_LCA})-(\ref{eq:gammaEOM_LCA}) analytically. During the
kinetic propagation $[\hat{H}_{\text{eff}}=T(\hat{p})$, i.e., $V=0]$,
Eqs.~(\ref{eq:qEOM_LCA})-(\ref{eq:gammaEOM_LCA}) have the analytical solution
\begin{align}
q_{t}  &  =q_{0}+t\,m^{-1}\cdot p_{0},\label{qEOM_kin_prop}\\
p_{t}  &  =p_{0},\label{pEOM_kin_prop}\\
A_{t}  &  =(A_{0}^{-1}+t\,m^{-1})^{-1},\label{eq:f_T_AEOM_LCA}\\
\gamma_{t}  &  =\gamma_{0}+t\,T(p_{0})\nonumber\\
&  \quad+(i\hbar/2)\,\text{ln}\,[\,\text{det}\,(I_{D}+t\,m^{-1}\cdot A_{0})],
\label{eq:f_T_gammaEOM_LCA}%
\end{align}
where $I_{D}$ is the $D\times D$ identity matrix, and during the potential
propagation $[\hat{H}_{\text{eff}}=V_{\text{eff}}(\hat{q})$, i.e., $m^{-1}%
=0]$, Eqs.~(\ref{eq:qEOM_LCA})-(\ref{eq:gammaEOM_LCA}) have the analytical
solution
\begin{align}
q_{t}  &  =q_{0},\\
p_{t,i}  &  =p_{0,i}-t\,V^{\prime}(q_{0})_{i}-t\,V^{\prime\prime\prime}%
(q_{0})_{ijk}\,\Sigma_{0,jk}/2,\\
A_{t}  &  =A_{0}-t\,V^{\prime\prime}(q_{0}),\\
\gamma_{t}  &  =\gamma_{0}-t\,V(q_{0}). \label{eq:VstepSO}%
\end{align}

By sequentially performing potential propagation for time $\Delta t/2$,
kinetic propagation for time $\Delta t$, and potential propagation for time
$\Delta t/2$, one obtains the second-order \textquotedblleft
potential-kinetic-potential\textquotedblright(VTV) algorithm. Swapping the
kinetic and potential propagation steps gives the second-order
\textquotedblleft kinetic-potential-kinetic\textquotedblright(TVT) algorithm.
Each of the two numerical algorithms gives the state $\psi_{t+\Delta t}$ at
time $t+\Delta t$ from the state $\psi_{t}$ at time $t$:
\begin{equation}
|\psi_{t+\Delta t}\rangle=\hat{U}_{2}(t+\Delta t,t;\psi)|\psi_{t}\rangle,
\label{eq:U2}%
\end{equation}
where $\hat{U}_{2}$ is the approximate second-order evolution operator
associated with the VTV or TVT algorithm. Note that these algorithms are
analogues of Faou and Lubich's second-order geometric integrator for the
variational GWD.~\cite{Faou_Lubich:2006, book_Hairer_Wanner:2006}

\subsection{High-order geometric integrators}

\label{sec-Compositions} An appropriate composition of the second order (VTV
or TVT) algorithm~(\ref{eq:U2}) yields integrators of higher order of
accuracy.~\cite{book_Hairer_Wanner:2006, book_Leimkuhler_Reich:2004,
Yoshida:1990, Suzuki:1990, McLachlan:1995, Wehrle_Vanicek:2011,
Sofroniou_Spaletta:2005, Choi_Vanicek:2019, Roulet_Vanicek:2019} More
precisely, any symmetric evolution operator $\hat{U}_{p}$ of even order $p$
yields an evolution operator $\hat{U}_{p+2}$ of order $p+2$ if it is
symmetrically composed as~\cite{book_Hairer_Wanner:2006}
\begin{align}
\hat{U}_{p+2}(t+\Delta t,t;\psi)  &  :=\hat{U}_{p}(t+\xi_{M}\Delta
t,t+\xi_{M-1}\Delta t;\psi)\nonumber\\
&  \qquad\qquad\cdots\hat{U}_{p}(t+\xi_{1}\Delta t,t;\psi),
\end{align}
where $M$ is the total number of composition steps and $\xi_{n}:=\sum
_{j=1}^{n}\gamma_{j}$ denotes the sum of the first $n$ real composition
coefficients $\gamma_{j}$. The choice of the composition coefficients must
satisfy the relations $\sum_{j=1}^{M}\gamma_{j}=1$ (consistency),
$\gamma_{M+1-j}=\gamma_{j}$ (symmetry), and $\sum_{j=1}^{M}\gamma_{j}^{p+1}=0$
(order increase guarantee).~\cite{book_Hairer_Wanner:2006} The simplest
symmetric composition methods are the recursive
triple-jump~\cite{Yoshida:1990} with $M=3$ and Suzuki's
fractal~\cite{Suzuki:1990} with $M=5$. Both methods can produce high-order
integrators. Because the number of composition steps and thus the
computational cost grows exponentially with the order of convergence, here we
only use the more efficient, non-recursive composition schemes that are, in a
certain sense, optimal. Kahan and Li~\cite{Kahan_Li:1997} obtained optimal
composition methods for the sixth and eighth orders by minimizing $\max
_{j}|\gamma_{j}|$, and Sofroniou and Spaletta~\cite{Sofroniou_Spaletta:2005}
found the tenth order optimal composition scheme by minimizing $\sum_{j=1}%
^{M}|\gamma_{j}|$. Suzuki's fractal gives the optimal fourth-order
scheme.~\cite{Choi_Vanicek:2019} For more details on the general as well as
optimal composition methods, see Ref.~\onlinecite{Choi_Vanicek:2019}.



\subsection{Geometric properties of the geometric integrators}

In the splitting algorithm described by Eqs.~(\ref{qEOM_kin_prop}%
)-(\ref{eq:VstepSO}), each potential or kinetic step of the local cubic
variational GWD is solved exactly. Therefore, before the composition, each
kinetic or potential step has all the geometric properties of the local cubic
variational method. All integrators, obtained by symmetric compositions of the
exact solutions~(\ref{qEOM_kin_prop})-(\ref{eq:VstepSO}) of the kinetic and
potential propagation steps, are time-reversible and conserve the norm and
symplectic structure (see Appendix~\ref{sec:num_symplectic_structure} for an
explicit proof of symplecticity). However, due to the splitting, they conserve
the effective energy only approximately, with an error $O(\Delta t^{M})$ where
the power $M$ is equal to or greater than the order of the algorithm. For
general proofs of these statements see Refs.~\onlinecite{book_Leimkuhler_Reich:2004,book_Hairer_Wanner:2006, Roulet_Vanicek:2019, Choi_Vanicek:2019, Roulet_Vanicek:2021a}.


\section{Numerical examples \label{sec:numerical_example}}

In the following, we investigate the local cubic variational GWD and the
proposed high-order integrators in several anharmonic model systems. Our
calculations show (i) the superior accuracy of the local cubic variational
method over the local harmonic GWD, (ii) the conservation of the effective
energy by the local cubic variational method, (iii) the preservation of the
geometric properties of the local cubic variational method by the geometric
integrators, and (iv) the efficiency of the high-order integrators. All values
are presented in natural units (n.u.): $\hbar=m=1$. We do not show the results
of the VTV algorithm and its compositions since they are almost identical to
the corresponding results of the TVT algorithm.

To analyze the error due to the local cubic approximation of the potential, we
compare the local cubic variational with the fully variational results. For
numerical experiments, we have therefore chosen the coupled Morse potential,
for which the expectation values of the potential energy, gradient, and
Hessian, needed in the fully variational GWD, can be computed analytically.

The nonseparable, $D$-dimensional \textquotedblleft coupled Morse
potential\textquotedblright~\cite{Fereidani_Vanicek:2023}
\begin{equation}
V(q)=V_{\text{eq}}+\sum_{j=1}^{D}V_{j}(q_{j})+V_{\text{cpl}}%
(q)\label{eq:Coupled-Morse}%
\end{equation}
consists of $D$ independent one-dimensional Morse potentials $V_{j}(q_{j})$,
which are coupled by a \textquotedblleft multidimensional Morse
coupling\textquotedblright\ $V_{\text{cpl}}(q)$, and $V_{\text{eq}}$ is the
potential at the equilibrium bond distance $q_{\text{eq}}$. Each
one-dimensional Morse potential
\begin{equation}
V_{j}(q_{j}):=d_{e}^{\prime}\,\big[1-y_{j}(a_{j}^{\prime},q_{j})\big]^{2}%
\label{eq:1D_Morse}%
\end{equation}
depends on the dissociation energy $d_{e}^{\prime}$ and one-dimensional Morse
variable~\cite{Braams_Bowman:2009}
\begin{equation}
y_{j}(a_{j}^{\prime},q_{j}):=\exp\big[-a_{j}^{\prime}\,(q_{j}-q_{\text{eq}%
,j})\big],\label{eq:Morse_var_1D}%
\end{equation}
which, in turn, is a function of the decay parameter $a_{j}^{\prime}$ and
$q_{j}$. The $D$-dimensional coupling term
\begin{equation}
V_{\text{cpl}}(q):=d_{e}\,\big[1-y(a,q)\big]^{2}\label{eq:cpl_Morse}%
\end{equation}
depends on the dissociation energy $d_{e}$ and $D$-dimensional Morse variable
\begin{equation}
y(a,q):=\exp\big[-a^{T}\cdot(q-q_{\text{eq}})\big],\label{eq:Morse_var}%
\end{equation}
which, in turn, is a function of the decay vector $a$ and $q$. The decay
parameter $a_{j}^{\prime}$, dissociation energy $d_{e}^{\prime}$, and
dimensionless anharmonicity $\chi_{j}^{\prime}$ are connected by the relation
$a_{j}^{\prime}=\chi_{j}^{\prime}\,\sqrt{8\,d_{e}^{\prime}}$%
.~\cite{Begusic_Vanicek:2019} Likewise, the decay vector $a$, dissociation
energy $d_{e}$, and dimensionless anharmonicity vector $\chi$ are related via
$a=\chi\,\sqrt{8\,d_{e}}$. The derivatives and expectation values of the
coupled Morse potential~(\ref{eq:Coupled-Morse}) are given in Ref.~\onlinecite{Fereidani_Vanicek:2023}.



\subsection{One-dimensional Morse potential \label{subsec:1D_Morse}}

\begin{figure}
\includegraphics[width=0.45\textwidth]{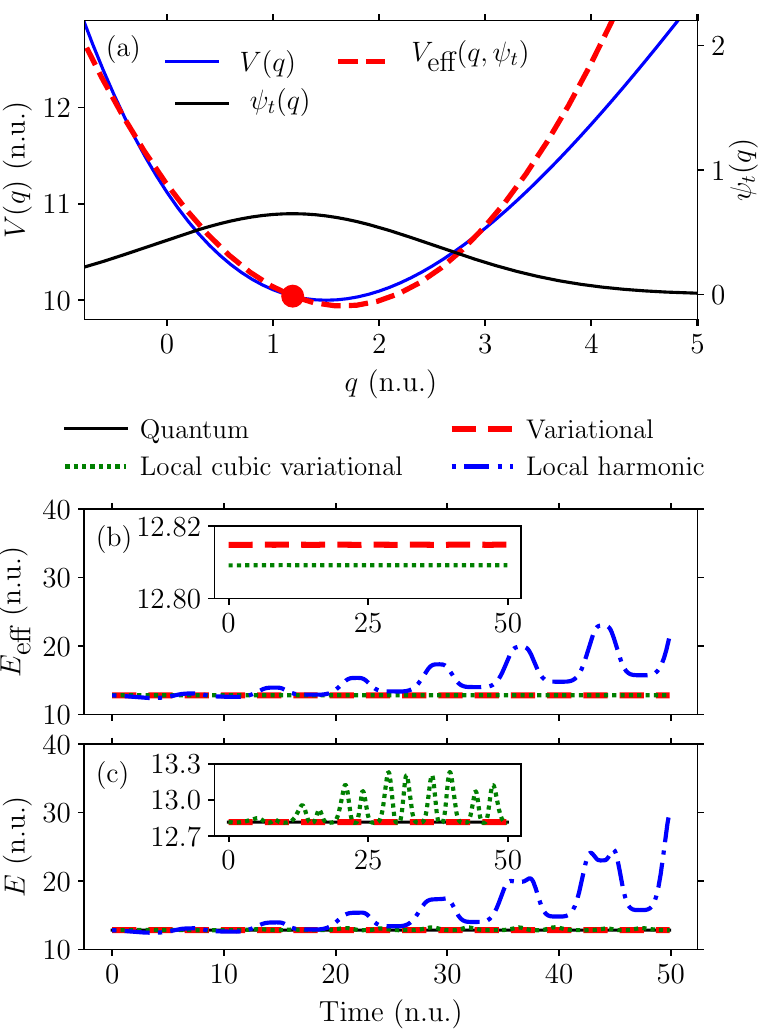}
\caption{Dynamics of a wavepacket $\psi_{t}$ in a one-dimensional Morse potential $V$.
Panel (a) shows the effective potential $V_{\textrm{eff}}$ of the local cubic variational GWD.
Panels (b) and (c) compare the effective $(E_{\textrm{eff}})$ and exact $(E)$ energies of the wavepacket
propagated with different methods. For clarity, the insets of panels (b) and (c) do not show the local harmonic results.  }\label{fig:1D_Morse}%

\end{figure}

\begin{figure}
\includegraphics[width=0.45\textwidth]{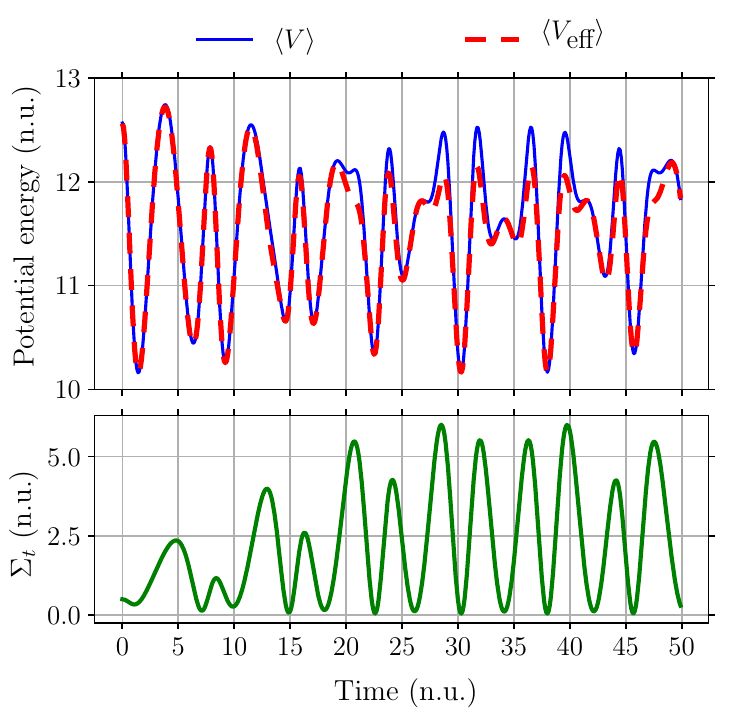}
\caption{Top panel compares the exact and effective potential energies of a  Gaussian wavepacket propagated
in a one-dimensional Morse potential with the local cubic variational GWD. These two energies are almost equal
whenever the position covariance, shown in the bottom panel, is small.}\label{fig:covq_LCA}%

\end{figure}

\begin{figure}
\includegraphics[width=0.45\textwidth]{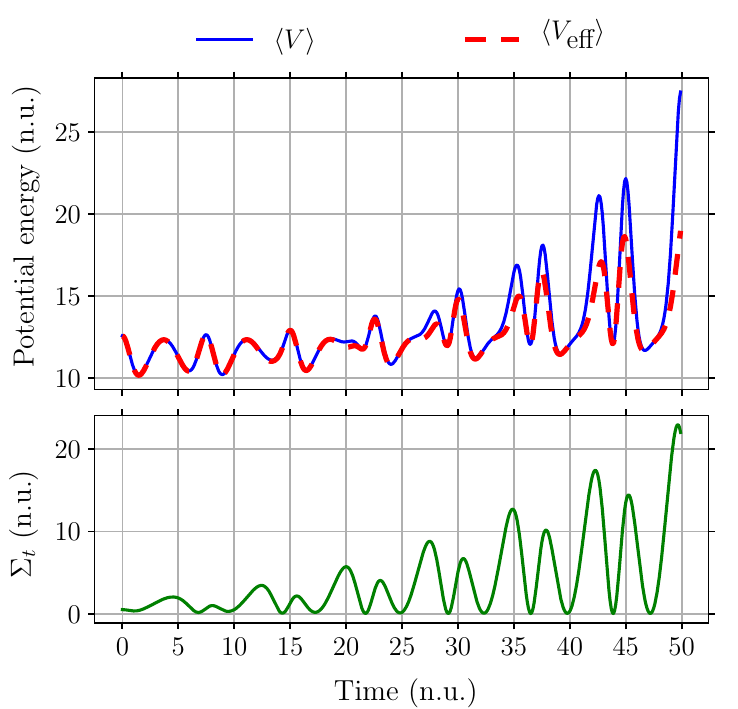}
\caption{Top panel compares the exact and effective potential energies of a Gaussian wavepacket
propagated in a one-dimensional Morse potential with the local harmonic GWD. These two energies are almost
identical whenever the position covariance, shown in the bottom panel, is small.}\label{fig:covq_TGA}%

\end{figure}

\begin{figure}
\includegraphics[width=0.45\textwidth]{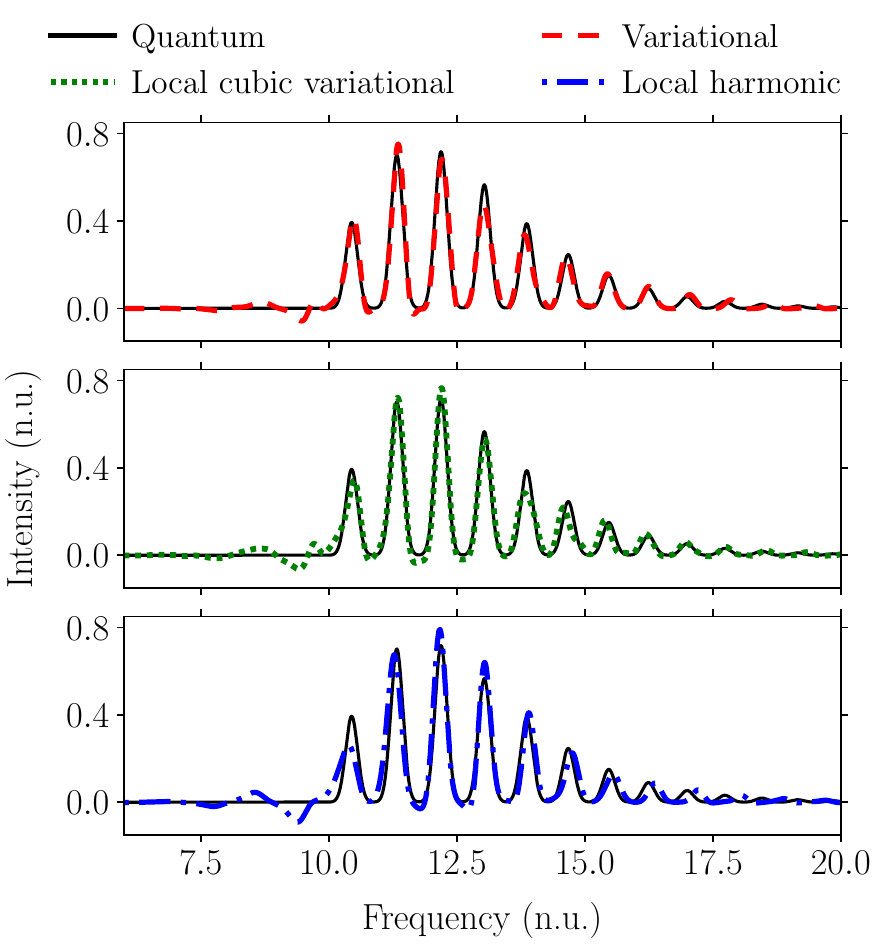}
\caption{Absorption spectrum of a one-dimensional Morse potential with anharmonicity $\chi=0.01$. The spectra obtained with the (a) variational, (b) local cubic variational, and (c) local harmonic GWD are compared to the exact quantum spectrum. }\label{fig:spec_1D}%

\end{figure}

Figures~\ref{fig:1D_Morse}-\ref{fig:spec_1D} analyze the dynamics of a
Gaussian wavepacket in a standard one-dimensional Morse potential, which is
also a special case of the coupled Morse potential~(\ref{eq:Coupled-Morse})
with zero coupling term ($d_{e}=0$) in Eq.~(\ref{eq:cpl_Morse}). The potential
parameters are $q_{\text{eq}}=1.5$, $V_{\text{eq}}=10$, $d_{e}^{\prime}=22.5$,
and anharmonicity $\chi^{\prime}=0.01$. The initial wavepacket was a real
Gaussian with $q_{0}=-0.5$, $p_{0}=0$, and $A_{0}=i\omega_{0}$, which is also
the ground vibrational state of a harmonic potential with frequency
$\omega_{0}=1$. This wavepacket was propagated for $5000$ steps of $\Delta
t=0.01$ with different GWD methods using the second-order geometric integrator
from Sec.~\ref{sec:geometric_integrators} or with the quantum dynamics using
the analogous second-order split-operator algorithm. The position grid for the
quantum dynamics consisted of $513$ equidistant ($\Delta q=30/512$) points
between $-5$ and $25$.

Panel (a) of Fig.~\ref{fig:1D_Morse} shows a snapshot at time $t=8.8$ of the
wavepacket propagated with the local cubic variational GWD. Note that the
associated effective potential $V_{\text{eff}}$ is not tangent to the exact
potential $V$, and thus differs from the local harmonic approximation. Panels
(b) and (c) compare the effective and exact energies of the wavepacket
propagated with different methods. The effective and exact energies of the
variational GWD are equal and conserved. Panel (b) shows that, unlike the
local harmonic method, the local cubic approach conserves its effective
energy. Panel (c) indicates that the energy of the local cubic method is not
conserved. However, its deviation from the exact value is significantly
smaller than that of the local harmonic one.

The expectation values of the exact and effective potential energies are
compared in Fig.~\ref{fig:covq_LCA} for the local cubic variational GWD and in
Fig.~\ref{fig:covq_TGA} for the local harmonic GWD. Both figures confirm that
the difference between the exact and effective potential energies increases
with time. However, this increase is more pronounced in the local harmonic
method. The lower panels of Figs.~\ref{fig:covq_LCA} and~\ref{fig:covq_TGA}
show the position covariance~(\ref{eq:cov_q}), which is basically the squared
width of a one-dimensional Gaussian wavepacket. In both figures, the exact and
effective potential energies are almost identical at local minima of the
position covariance, since at these times the Gaussian wavepacket becomes
infinitesimally narrow, and the potential is thus approximated over a smaller
region of space. In contrast, the exact and effective potential energies
differ the most at local maxima of the position covariance.

We also calculated the absorption spectrum, corresponding to the electronic
transition from a harmonic PES to the just described Morse PES, as the Fourier
transform~\cite{Heller:1978,Heller:1981a,book_Mukamel:1999, book_Tannor:2007,
Lami_Santoro:2004}
\begin{equation}
\sigma(\omega)=\text{Re}\int_{0}^{\infty}C(t)\,e^{i(\omega+E_{1,g}/\hbar
)t}\,dt \label{eq:spectra}%
\end{equation}
of the wavepacket autocorrelation function
\begin{equation}
C(t)=\langle\psi_{0}|\psi_{t}\rangle, \label{eq:auto_corr}%
\end{equation}
where $E_{1,g}$ is the vibrational zero point energy of state $\psi_{0}$ in
the harmonic PES before photon absorption. The overlap~(\ref{eq:auto_corr}) of
Gaussian wavepackets $\psi_{0}$ and $\psi_{t}$ is given by Eq. (29) in
Ref.~\onlinecite{Begusic_Vanicek:2020a} or by Eq. (B4) in Ref.~\onlinecite{Fereidani_Vanicek:2023}.

Figure~\ref{fig:spec_1D} compares spectra computed with different methods. To
obtain smooth spectra, the simulation time in Fig.~\ref{fig:spec_1D} was ten
times longer than the simulation time in Fig.~\ref{fig:1D_Morse}. A Gaussian
damping function with a half-width at half-maximum of $10$ was applied to each
autocorrelation function before computing the spectrum by Fourier transform.
All three semiclassical methods perform well in the calculation of spectra,
but the nonlinear character of their effective Hamiltonian~(\ref{eq:Heff})
leads to small, yet visible unphysical negative intensities. As expected, the
fully variational spectrum is slightly more accurate than the local cubic
variational one, and its negative features are smaller. Similarly, the local
cubic variational spectrum is more accurate than the local harmonic one,
especially in the low-frequency region and in the high-frequency tail.


\subsection{Two-dimensional coupled Morse potential
\label{subsec:2D_coupled_Morse}}

\begin{figure}
\includegraphics[width=0.45\textwidth]{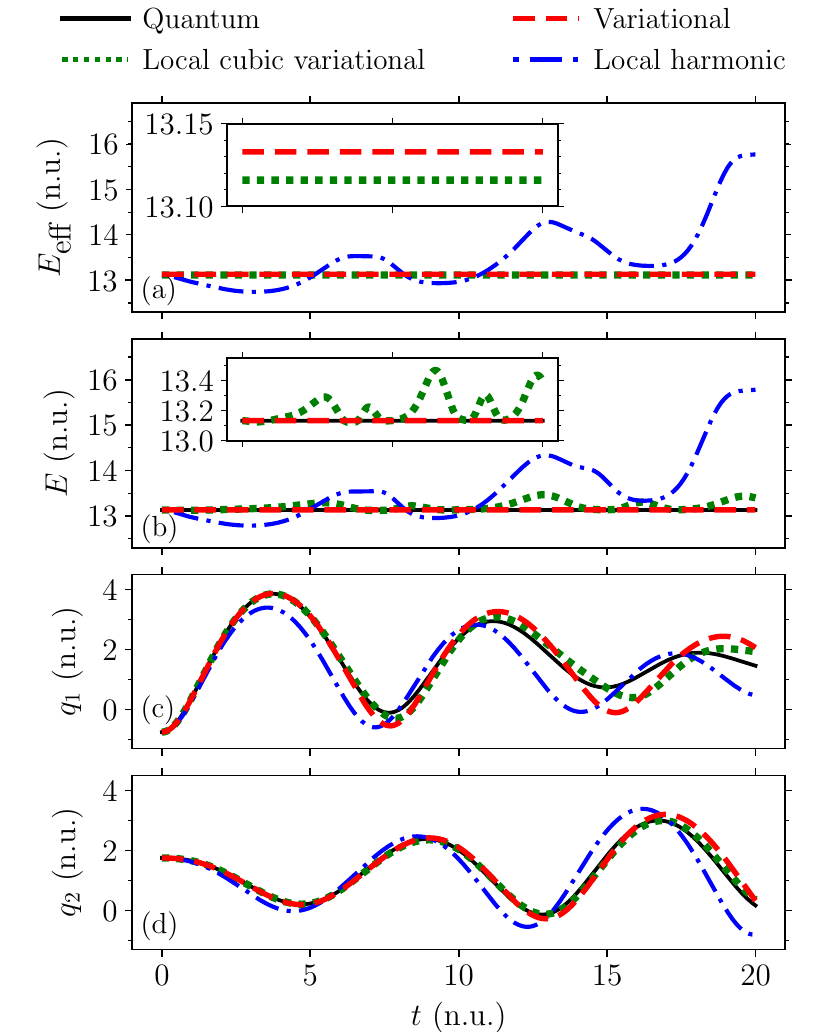}
\caption{Dynamics of a wavepacket propagated in a two-dimensional coupled Morse potential with different methods. Effective energy [panel (a)],  energy [panel (b)], and position of the Gaussian's center along two different coordinates [panels (c) and (d)] are shown. For clarity, the insets of panels (a) and (b) do not show the local harmonic results.}\label{fig:2D_coupled_Morse}%

\end{figure}
\begin{figure}
\includegraphics[width=0.45\textwidth]{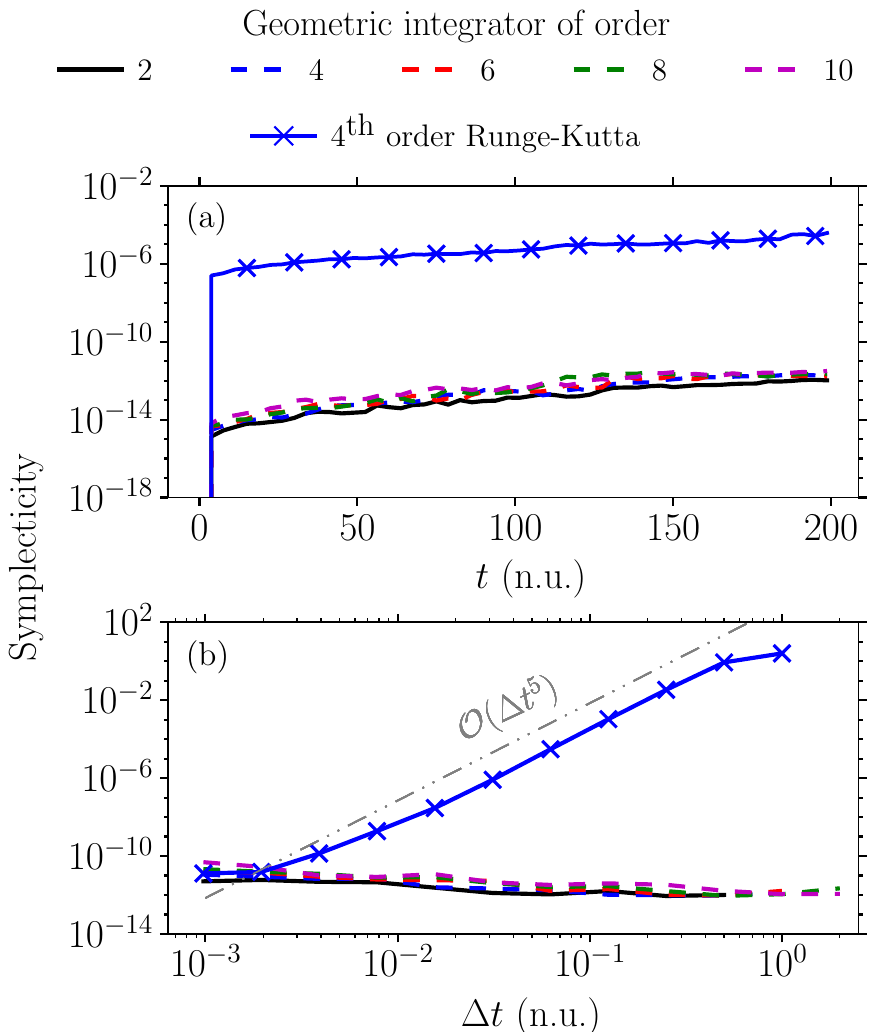}
\caption{Conservation of the symplectic structure of the Gaussian wavepackets by the geometric integrators and
its nonconservation by  the fourth-order Runge-Kutta method. The system is the same two-dimensional coupled
Morse potential as in Fig.~\ref{fig:2D_coupled_Morse}. The symplecticity~(\ref{eq:symplecticity}) is shown (a) as
a function of time $t$ for a fixed and fairly large time step $\Delta t= 2^{-4}\, \textrm{n.u.}$ and (b) as a function
of time step $\Delta t$ at the final time $t_{f}=200\, \textrm{n.u.}$
}\label{fig:symplecticity}
\end{figure}

Figures~\ref{fig:2D_coupled_Morse} and~\ref{fig:symplecticity} analyze the
dynamics of a wavepacket in a two-dimensional coupled Morse
potential~(\ref{eq:Coupled-Morse}) with energy $V_{\text{eq}}=10$ at the
equilibrium position $q_{\text{eq}}=(1,1)$. The potential is composed of two
one-dimensional Morse potentials with the same dissociation energy
$d_{e}^{\prime}=11.25$ and different anharmonicities $\chi_{1}^{\prime}=0.02$
and $\chi_{2}^{\prime}=0.017$. The parameters of the coupling term are
$d_{e}=5.75$ and $\chi=(0.014,0.017)$. The initial wavepacket was a real
two-dimensional Gaussian~(\ref{eq:GWP}) with $q_{0}=(-0.75,1.75)$,
$p_{0}=(0,0)$, and a diagonal width matrix $A_{0}$ with non-zero elements
$A_{0,11}=A_{0,22}=i$. This wavepacket was propagated for $20000$ steps of
$\Delta t=0.001$ with different methods using the second-order geometric or
split-operator integrator. The position grid for the exact quantum dynamics
consisted of $257$ equidistant $(\Delta q=20/256)$ points between $-5$ and
$15$ in both directions.

Panels (a) and (b) of Fig.~\ref{fig:2D_coupled_Morse} confirm our observations
for the one-dimensional Morse potential in Fig.~\ref{fig:1D_Morse}: (i) the
effective and exact energies of the fully variational method are conserved and
equal to each other and to the quantum energy, (ii) the local cubic
variational GWD conserves the effective energy, but not the exact energy, and
(iii) the local harmonic GWD conserves neither the effective nor the exact
energy. Panels (c) and (d) show the expectation values of the position of the
wavepacket. For very short times, all approximations give accurate quantum
results, but their accuracy decreases with time. However, the variational and
local cubic variational methods remain accurate longer than the local harmonic
GWD. The spectra of this two-dimensional coupled Morse potential, calculated
with different methods, are available in the supplementary material.



In Fig.~\ref{fig:symplecticity}, we analyze the symplecticity of the geometric
integrators for the local cubic variational GWD. We verified the symplecticity
numerically by measuring the distance
\begin{align}
d  &  \big(\Phi_{t}^{\prime}(z_{0})^{T}\cdot\omega(z_{t})\cdot\Phi_{t}%
^{\prime}(z_{0}),\omega(z_{0})\big)\nonumber\\
&  \equiv\lVert\Phi_{t}^{\prime}(z_{0})^{T}\cdot\omega(z_{t})\cdot\Phi
_{t}^{\prime}(z_{0})-\omega(z_{0})\rVert\label{eq:symplecticity}%
\end{align}
between the \textquotedblleft initial" and \textquotedblleft final" symplectic
structure matrices $\omega(z_{0})$ and $\omega(z_{t})$. Here, the vector $z$
contains elements of the Gaussian's parameters propagated by a numerical
integrator with flow $\Phi_{t}(z)$ and Jacobian $\Phi_{t}^{\prime}(z)$, and
$\omega(z)$ is a skew-symmetric matrix representing the non-canonical
symplectic two-form of Gaussian wavepackets. Hagedorn proposed a different but
equivalent parametrization of the Gaussian wavepacket,~\cite{Hagedorn:1980_v2}
in which the equations of motion of the GWD become simpler (see
Appendix~\ref{sec:Hagedorn_rep}). We, therefore, evaluated the Jacobian and
symplecticity in Hagedorn's parametrization; see
Appendix~\ref{sec:num_symplectic_structure} for more details.

Figure~\ref{fig:symplecticity} shows the symplecticity~(\ref{eq:symplecticity}%
) of the Gaussian wavepacket propagated with the local cubic variational GWD
in the two-dimensional coupled Morse potential analyzed in
Fig.~\ref{fig:2D_coupled_Morse}. Although the simulation time in
Fig.~~\ref{fig:symplecticity} is ten times longer than the simulation time in
Fig.~\ref{fig:2D_coupled_Morse}, all geometric integrators conserve the
non-canonical symplectic structure exactly, regardless of the size of the time
step. In contrast, the figure clearly shows that the fourth-order Runge-Kutta
method~\cite{book_Leimkuhler_Reich:2004} does not conserve the symplectic
structure of the local cubic variational GWD.


\subsection{Twenty-dimensional coupled Morse potential
\label{subsec:20D_coupled_Morse}}

\begin{figure}
\includegraphics[width=0.45\textwidth]{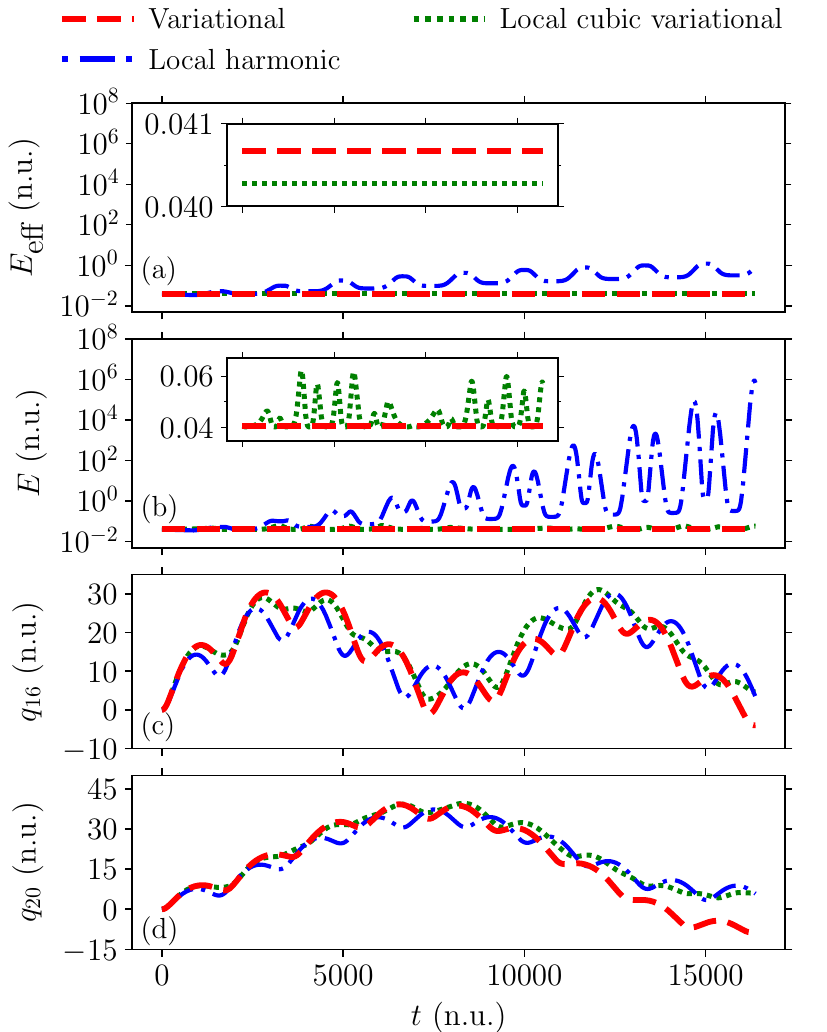}
\caption{Dynamics of a Gaussian wavepacket propagated in a twenty-dimensional coupled Morse potential with various GWD methods. Effective energy [panel (a)],  energy [panel (b)], and position of the Gaussian's center along two different coordinates [panels (c) and (d)] are shown. For clarity, the insets of panels (a) and (b) do not show the local harmonic results.}\label{fig:20D_Coupled_Morse}%

\end{figure}

\begin{figure}
\includegraphics[width=0.45\textwidth]{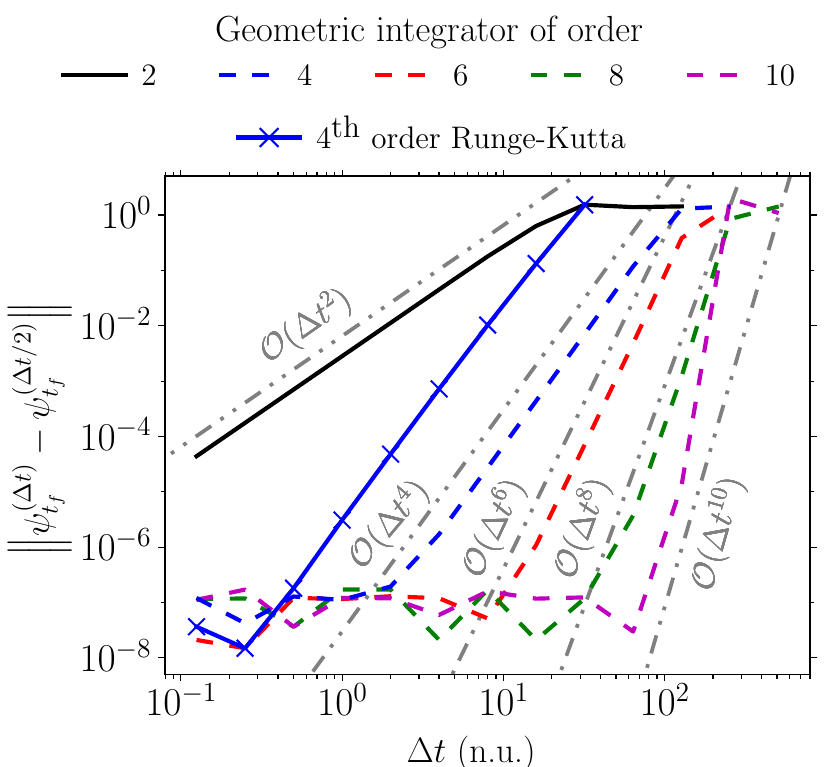}
\caption{Convergence of various integrators for the local cubic variational GWD, measured by the convergence error at the final time $t_{f} = 2^{16} \, \textrm{n.u.}=65536 \, \textrm{n.u.}$ as a function of the time step $\Delta t$.
The convergence error is defined as the distance $d(\psi^{(\Delta t)}_{t},\psi^{(\Delta t/2)}_{t})\equiv \lVert \psi^{(\Delta t)}_{t}-\psi^{(\Delta t/2)}_{t} \rVert$, where $\psi^{(\Delta t)}_{t}$ denotes the state at time $t$ obtained after propagation with time step $\Delta t$.
}\label{fig:convergence}
\end{figure}

\begin{figure}
\includegraphics[width=0.45\textwidth]{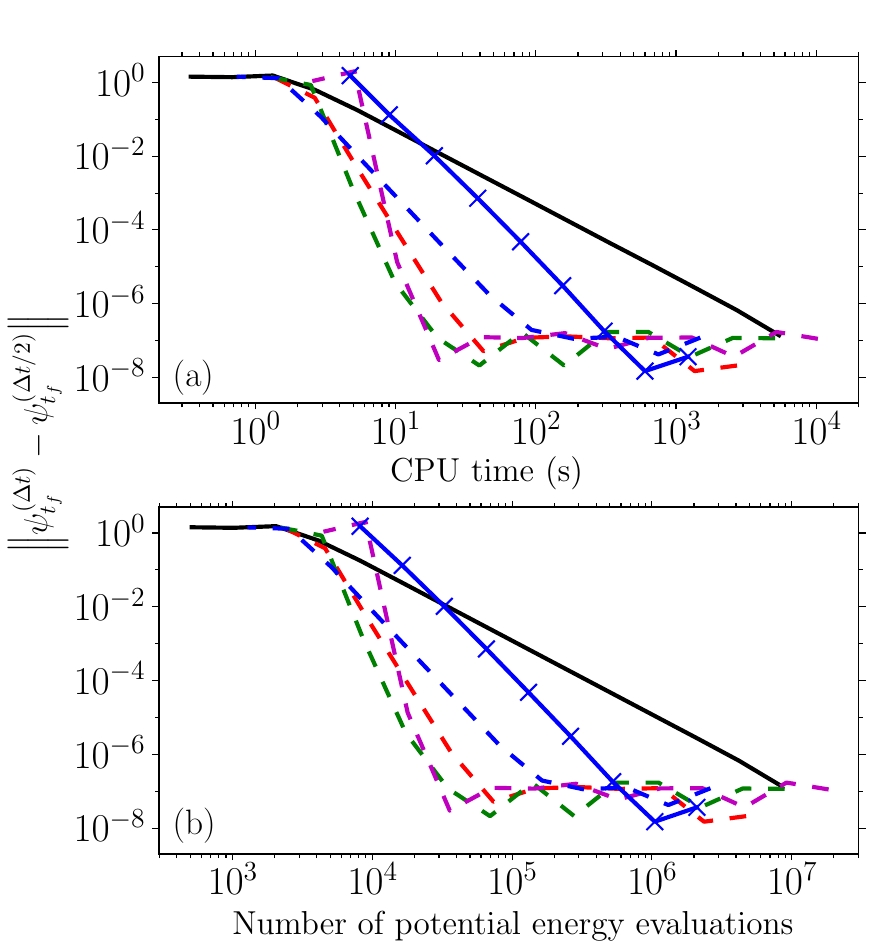}
\caption{ Efficiency of various integrators  for the local cubic variational GWD. Efficiency is measured
by plotting the convergence error, defined in the caption of Fig. \ref{fig:convergence}, as a function
of (a) the computational cost (CPU time) or (b) the number of potential energy evaluations.
Line labels are the same as those in Fig. \ref{fig:convergence}.
}\label{fig:efficiency}

\end{figure}

\begin{figure}
\includegraphics[width=0.45\textwidth]{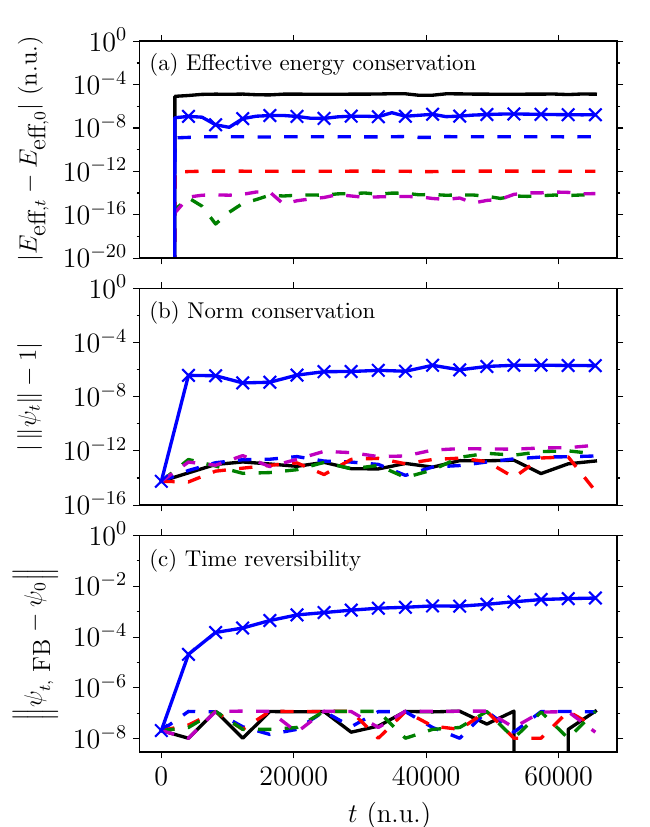}
\caption{ Geometric properties of various integrators for the local cubic variational GWD as a function of time $t$ for a relatively large time step $\Delta t = 8 \, \textrm{n.u.}$ (a) Effective energy, (b) norm, and (c) time reversibility [Eq.~(\ref{eq:TR})] are shown. Line labels are the same as those in Fig. \ref{fig:convergence}.
}\label{fig:geom_prop_t}
\end{figure}

\begin{figure}
\includegraphics[width=0.45\textwidth]{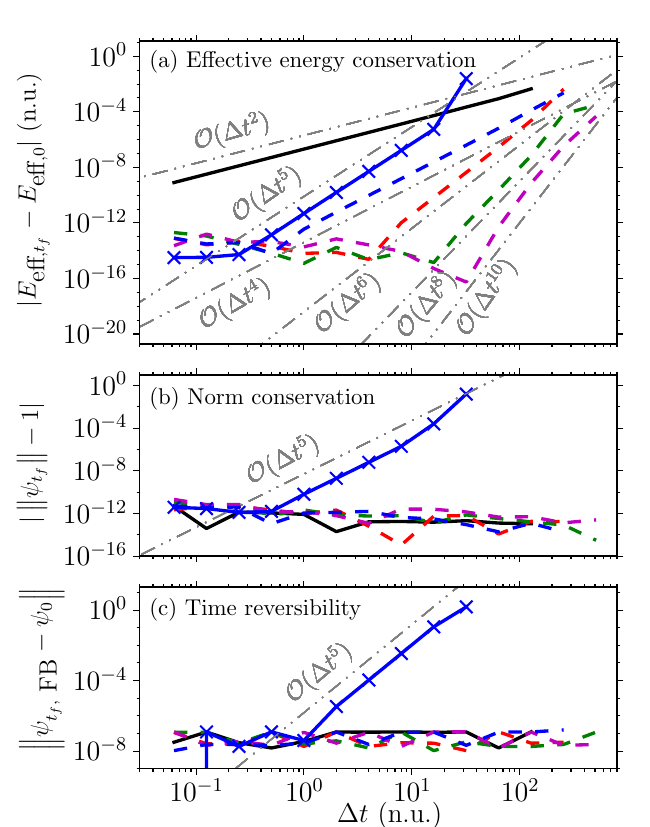}
\caption{Geometric properties of various integrators for the local cubic variational GWD  as a function of the time step $\Delta t$ measured at the final time $t_{f} = 2^{16} \, \textrm{n.u.}=65536 \, \textrm{n.u.}$ (a) Effective energy, (b) norm, and (c) time reversibility [Eq.~(\ref{eq:TR})] are shown. Line labels are the same as those in Fig. \ref{fig:convergence}.
}\label{fig:geom_prop_dt}
\end{figure}

Figures~\ref{fig:20D_Coupled_Morse}-\ref{fig:geom_prop_dt} analyze the
dynamics of a Gaussian wavepacket in a twenty-dimensional coupled Morse
potential~(\ref{eq:Coupled-Morse}) composed of twenty one-dimensional Morse
potentials~(\ref{eq:1D_Morse}) with the same dissociation energy
$d_{e}^{\prime}=0.1$ and with anharmonicity parameters $\chi_{j}^{\prime
}\,\,(j=1,\dots,20)$ uniformly varying in the range between $0.001$ and
$0.005$ $(\Delta\chi^{\prime}=0.004/19)$. Parameters of the coupling
term~(\ref{eq:cpl_Morse}) were $d_{e}=0.075$ and $\chi_{j}=(3/4)\,\chi
_{j}^{\prime}$, and the energy at the equilibrium position $q_{\text{eq}%
,j}=10\,\,(j=1,\dots,20)$ was $V_{\text{eq}}=0$. The twenty-dimensional real
initial Gaussian wavepacket had zero position, zero momentum, and a diagonal
width matrix with non-zero elements $A_{jj}=4\,d_{e}\chi_{j}\,i$.

In Fig.~\ref{fig:20D_Coupled_Morse}, this wavepacket was propagated for
$2^{17}=131072$ steps of $\Delta t=0.125$ with various GWD methods using the
second-order geometric integrator. Unlike Fig.~\ref{fig:2D_coupled_Morse},
Fig.~\ref{fig:20D_Coupled_Morse} lacks the exact quantum benchmark because the
system is not solvable by the grid-based methods due to its high
dimensionality. Even using the multi-configurational time-dependent Hartree
method,~\cite{Meyer_Cederbaum:1990,Beck_Jackle:2000} exact quantum calculation
for this high-dimensional system would be very difficult. Panels (a) and (b)
of Fig.~\ref{fig:20D_Coupled_Morse} show the exact and effective energies of
the wavepacket propagated with different semiclassical methods. The effective
and exact energies of the fully variational method are equal and conserved.
Both the effective and exact energies of the wavepacket propagated with the
local cubic variational or local harmonic GWD differ from the corresponding
energies obtained by the variational GWD. However, the difference is much
greater for the local harmonic GWD, especially at longer times. Furthermore,
unlike the local harmonic GWD, the local cubic variational GWD conserves its
effective energy. Panels (c) and (d) display the position of the Gaussian's
center. For very short times, all semiclassical methods overlap almost
perfectly. However, the local cubic variational results remain close to the
fully variational results for longer times.


To analyze the convergence and geometric properties of the integrators, we
repeated the local cubic variational simulation with several high-order
geometric integrators and with the fourth-order Runge-Kutta method.

Convergence rates of various methods are compared in
Fig.~\ref{fig:convergence}. For all methods, the obtained orders of
convergence agree with the predicted ones, indicated by the gray straight
lines. The plateau indicates the machine precision error. Since the high-order
methods require a large number of composition substeps to be carried out at
each time step $\Delta t$, a higher order of convergence does not guarantee
higher efficiency. Figure~\ref{fig:efficiency} measures the efficiency
directly by plotting the convergence error as a function of either the CPU
time or the number of potential energy evaluations. The similarity between
panels (a) and (b) confirms that the cost of initialization and finalization
is negligible to the cost of potential propagation; this will be even more the
case in expensive on-the-fly ab initio applications. In addition,
Fig.~\ref{fig:efficiency} shows that high-order geometric integrators are more
efficient than both the second-order geometric integrator and the fourth-order
Runge-Kutta method. For example, below a rather large error of $10^{-1}$, the
eighth-order integrator is more efficient than the sixth-, fourth-, and
second-order integrators, and for a moderate error of $10^{-4}$, the
fourth-order geometric integrator is almost $10$ times faster than the
second-order geometric algorithm and almost $3.5$ faster than the fourth-order
Runge-Kutta method. As a result, the use of high-order geometric integrators
can simultaneously increase both accuracy and efficiency.

Figure \ref{fig:geom_prop_t} shows how the effective energy, norm, and time
reversibility depend on time, whereas Fig.~\ref{fig:geom_prop_dt} analyzes how
these geometric properties depend on the time step. Panels~(a) of both figures
analyze the conservation of the effective energy measured by the difference
$|E_{\text{eff},t}-E_{\text{eff},0}|$. Although the local cubic variational
GWD conserves the effective energy, the splitting nature of the geometric
integrators reduces this exact conservation.~\cite{book_Hairer_Wanner:2006,
Roulet_Vanicek:2019} This is because the alternation between kinetic [$\hat
{H}_{\mathrm{eff}}=T(\hat{p})$] and potential [$\hat{H}_{\mathrm{eff}%
}=V_{\mathrm{eff}}(\hat{q})$] propagations makes the effective Hamiltonian
time-dependent. The geometric integrators and the fourth-order Runge-Kutta
method conserve the effective energy only approximately, with the order in
$\Delta t$ equal to or greater than their order of convergence. Conservation
of the norm is measured by the deviation $|\lVert\psi_{t}\rVert-1|$, where
\begin{equation}
\lVert\psi_{t}\rVert=\text{exp}(-\text{Im}\,\gamma_{t}/\hbar)\,\big[\text{det}%
(\pi\hbar/\text{Im}\,A_{t})\big]^{1/4} \label{eq:Norm}%
\end{equation}
is the norm of the Gaussian wavepacket~(\ref{eq:GWP}). Unlike the effective
energy, which is conserved only approximately, panels~(b) of
Figs.~\ref{fig:geom_prop_t} and~\ref{fig:geom_prop_dt} confirm the exact norm
conservation by the geometric integrators. The fourth-order Runge-Kutta
method, on the other hand, fails to conserve the norm unless it is fully
converged, which requires using a very small time step. In panels (c) of
Figs.~\ref{fig:geom_prop_t} and~\ref{fig:geom_prop_dt}, the time reversibility
is checked by plotting the distance
\begin{equation}
d\big(\psi_{t,\text{FB}},\psi_{0}\big)\equiv\lVert\psi_{t,\text{FB}}-\psi
_{0}\rVert\label{eq:TR}%
\end{equation}
between the initial state $\psi_{0}$ and the \textquotedblleft
forward-backward\textquotedblright\ propagated state $\psi_{t,\text{FB}}$
(i.e., the state propagated first forward in time for time $t$ and then
backward in time for time $t$). Unlike the fourth-order Runge-Kutta method,
the geometric integrators are exactly time-reversible. Note that the
convergence of norm, energy, and reversibility by the fourth-order Runge-Kutta
method appears to be $\mathcal{O}(\Delta t^{5})$, i.e., faster than
$\mathcal{O}(\Delta t^{4})$.

We do not repeat the analysis of symplecticity, which was analyzed for the
two-dimensional coupled Morse potential in Fig.~\ref{fig:symplecticity},
because its evaluation for the twenty-dimensional system is unnecessarily
expensive. Finally, as any exact solution of a genuinely nonlinear TDSE, the
local cubic variational GWD conserves neither the inner product nor the
distance between two states (see Fig.~S2 in the supplementary material).


\section{Conclusion \label{sec:conclusion}}

We have shown that applying the local cubic approximation to the potential in
the variational GWD improves the accuracy over Heller's original GWD, which is
based only on the local harmonic approximation. In contrast to the original
GWD, the local cubic variational method is symplectic and conserves the
effective energy. Most importantly, the local cubic method is also a practical
implementation of the variational method for future on-the-fly ab initio
applications, since it avoids the difficult evaluation of $\langle\hat
{V}\rangle$, $\langle\hat{V}^{\prime}\rangle$, and $\langle\hat{V}%
^{\prime\prime}\rangle$. Although more accurate than the local harmonic GWD,
the local cubic variational method significantly increases the computational
cost because it requires the third derivative of the potential along the
trajectory. To partially ameliorate this issue, we accelerated the method by
designing and implementing high-order geometric integrators. We showed that
these high-order integrators not only increased the efficiency and accuracy
over the second-order algorithm, but also preserved most geometric properties
of the local cubic variational GWD. Although the local cubic variational
method is only a crude approximation of the exact solution of the TDSE in
anharmonic systems, our well-converged results obtained with the high-order
integrators ensure that the numerical errors are negligible.

In conclusion, the local cubic variational GWD is promising for future
on-the-fly ab initio applications to medium-sized molecules. We expect it to
account more accurately for anharmonic effects, especially in ultrafast
spectroscopy, where only short-time dynamics is required. Its potential to
approximately include tunneling also deserves further attention.

\section*{Supplementary material}

See the supplementary material for the absorption spectra of the
two-dimensional coupled Morse potential from
Sec.~\ref{subsec:2D_coupled_Morse} and for the numerical demonstration of the
nonconservation of the inner product and distance by the local cubic
variational GWD.

\section*{Acknowledgments}

The authors acknowledge the financial support from the European Research
Council (ERC) under the European Union's Horizon 2020 research and innovation
programme (grant agreement No. 683069 -- MOLEQULE).

\section*{Data availability}

The data that support the findings of this study are available within the
article and its supplementary material.

\appendix

\section{\label{sec:Hagedorn_rep}Gaussian wavepacket dynamics in Hagedorn's
parametrization}

In Hagedorn's parametrization,~\cite{Hagedorn:1980_v2} the Gaussian wavepacket
(\ref{eq:GWP}) is written as
\begin{align}
&  \psi_{t}(q)= (\pi\hbar)^{-D/4} (\text{det} \, Q_{t})^{-1/2}\nonumber\\
&  \quad\times\text{exp}\bigg[\frac{i}{\hbar}\bigg(\frac{1}{2}\,x^{T} \cdot
P_{t} \cdot Q^{-1}_{t}\cdot x +p_{t}^{T} \cdot x+S_{t}\bigg)\bigg],
\label{eq:HGWP}%
\end{align}
where the new parameters $Q_{t}$ and $P_{t}$ are two $D \times D$ complex
matrices, related to the width of the Gaussian by $A_{t}=P_{t} \cdot
Q^{-1}_{t}$ and satisfying the relations~\cite{Hagedorn:1980_v2,
book_Lubich:2008}
\begin{align}
&  Q^{T}_{t} \cdot P_{t} - P^{T}_{t} \cdot Q_{t}=0,\\
&  Q^{\dagger}_{t} \cdot P_{t} - P^{\dagger}_{t} \cdot Q_{t}= 2 i I_{D},
\label{eq:QPrels}%
\end{align}
and $S_{t}$ is a real scalar that generalizes the classical action.

Equations of motion for $\dot{q}_{t}$ and $\dot{p}_{t}$ [Eqs.~(\ref{eq:qEOM})
and~(\ref{eq:pEOM})] remain unchanged, whereas the equations for $\dot{A}_{t}$
and $\dot{\gamma}_{t}$ [Eqs.~(\ref{eq:AEOM}) and~(\ref{eq:gammaEOM})] are
replaced by
\begin{align}
\dot{Q}_{t}  &  =m^{-1}\cdot P_{t},\label{eq:H_QEOM_LCA}\\
\dot{P}_{t}  &  =-V_{2}\cdot Q_{t},\label{eq:H_PEOM_LCA}\\
\dot{S}_{t}  &  =T(p_{t})-V_{0}. \label{eq:H_SEOM_LCA}%
\end{align}

\subsection{\label{subsec:H_geometric_integrators}Geometric integrators in
Hagedorn's parametrization}

In Hagedorn's parametrization, the kinetic propagation flow $\Phi_{\text{T}%
,t}$ is~\cite{Vanicek:2023}
\begin{align}
q_{t}  &  =q_{0}+t\,m^{-1}\cdot p_{0},\label{eq:Hf_T_qEOM_LCA}\\
p_{t}  &  =p_{0},\label{eq:Hf_T_pEOM_LCA}\\
Q_{t}  &  =Q_{0}+t\,m^{-1}\cdot P_{0},\label{eq:Hf_T_QEOM_LCA}\\
P_{t}  &  =P_{0},\label{eq:Hf_T_PEOM_LCA}\\
S_{t}  &  =S_{0}+t\,T(p_{0}), \label{eq:Hf_T_SEOM_LCA}%
\end{align}
whereas the potential propagation flow $\Phi_{\text{V},t}$
is~\cite{Vanicek:2023}
\begin{align}
q_{t}  &  =q_{0},\label{eq:Hf_V_qEOM_LCA}\\
p_{t}  &  =p_{0}-t\,V_{1}(q_{0},Q_{0}),\label{eq:Hf_V_pEOM_LCA}\\
Q_{t}  &  =Q_{0},\label{eq:Hf_V_QEOM_LCA}\\
P_{t}  &  =P_{0}-t\,V_{2}(q_{0},Q_{0})\cdot Q_{0},\label{eq:Hf_V_PEOM_LCA}\\
S_{t}  &  =S_{0}-t\,V_{0}(q_{0},Q_{0}). \label{eq:Hf_V_SEOM_LCA}%
\end{align}

\section{Symplecticity of the geometric integrators}


\label{sec:num_symplectic_structure}

To verify the symplecticity of a numerical integrator designed for the GWD,
one should show that the Jacobian $\Phi_{t}^{\prime}(z_{0})$ of the integrator
with the flow~$z_{t}=\Phi_{t}(z_{0})$ satisfies the
condition~\cite{book_Hairer_Wanner:2006}
\begin{equation}
\Phi_{t}^{\prime}(z_{0})^{T}\cdot\omega(z_{t})\cdot\Phi_{t}^{\prime}%
(z_{0})=\omega(z_{0}), \label{eq:JBJ=B}%
\end{equation}
where $z_{t}$ is a vector containing the Gaussian's parameters and
$\omega(z_{t})$ is the symplectic structure of the Gaussian wavepackets. [Note
that matrix $\omega$ in Eq.~(\ref{eq:JBJ=B}) is the inverse of matrix $B$
appearing in Eq. (4.2) of Chapter VII of Ref.~\onlinecite{book_Hairer_Wanner:2006}).]

In Hagedorn's parametrization, the reduced symplectic form
(\ref{eq:sympl_Agamma}) becomes the standard symplectic
form~\cite{Ohsawa:2015b}
\begin{align}
\omega(z)  &  =dq_{j}\wedge dp_{j}+(\hbar/2)\,dQ_{kl}^{(1)}\wedge
dP_{kl}^{(1)}\nonumber\\
&  \quad+(\hbar/2)\,dQ_{kl}^{(2)}\wedge dP_{kl}^{(2)} \label{eq:sympl_QPS}%
\end{align}
with
\[
z_{t}:=\bigg(q_{t}^{T},p_{t}^{T},\widetilde{Q_{t}^{(1)}}^{T},\widetilde{P_{t}%
^{(1)}}^{T},\widetilde{Q_{t}^{(2)}}^{T},\widetilde{P_{t}^{(2)}}^{T}%
\bigg)^{T}\in\mathbb{R}^{2D+4D^{2}},\label{eq:z_t}%
\]
where $Q_{t}^{(1)}$ and $Q_{t}^{(2)}$ are real and imaginary parts of matrix
$Q_{t}=Q_{t}^{(1)}+iQ_{t}^{(2)}$, and $P_{t}^{(1)}$ and $P_{t}^{(2)}$ are real
and imaginary parts of matrix $P_{t}=P_{t}^{(1)}+iP_{t}^{(2)}$. Moreover,
$\widetilde{\Lambda}$, which is used for $\Lambda=Q_{t}^{(1)},P_{t}%
^{(1)},Q_{t}^{(2)}$, and $P_{t}^{(2)}$, is a $D^{2}$-dimensional vector
containing elements of the $D\times D$ matrix $\Lambda$ in a column-by-column
manner, i.e., $\widetilde{\Lambda}_{j+D(k-1)}=\Lambda_{jk}$. Matrix
representation of (\ref{eq:sympl_QPS}) is the constant $(2D+4D^{2}%
)$-dimensional block-diagonal matrix
\begin{equation}
\omega(z_{t})=%
\begin{pmatrix}
J_{2D} & 0 & 0\\
0 & \frac{\hbar}{2}J_{2D^{2}} & 0\\
0 & 0 & \frac{\hbar}{2}J_{2D^{2}}%
\end{pmatrix}
. \label{eq:HOmega}%
\end{equation}
In Eq.~(\ref{eq:HOmega}), $J_{2D^{2}}=J_{2D}\otimes I_{D}$ and $J_{2D}%
=J_{2}\otimes I_{D}$, where
\begin{equation}
J_{2}=%
\begin{pmatrix}
0 & 1\\
-1 & 0
\end{pmatrix}
\label{eq:J}%
\end{equation}
is the canonical $2\times2$ symplectic matrix.

The flow of a geometric integrator of any order is composed of a sequence of
kinetic and potential flows. Its Jacobian is therefore equal to the matrix
product of the Jacobians of the elementary flows. To demonstrate the
symplecticity of geometric integrators, we first analyze the symplecticity of
the kinetic and potential flows.

\subsection{Symplecticity of the kinetic flow}

The kinetic flow $\Phi_{\text{T},t}(z_{0})$ given by
Eqs.~(\ref{eq:Hf_T_qEOM_LCA})-(\ref{eq:Hf_T_PEOM_LCA}) does not depend on the
coefficients $V_{0}$, $V_{1}$, and $V_{2}$, and thus is identical for all GWD
methods. The Jacobian of this flow is the block-diagonal
matrix~\cite{Fereidani_Vanicek:2023}
\begin{equation}
\Phi_{\text{T},t}^{\prime}(z_{0})=%
\begin{pmatrix}
M_{2D} & 0 & 0\\
0 & M_{2D^{2}} & 0\\
0 & 0 & M_{2D^{2}}%
\end{pmatrix}
, \label{eq:HJT}%
\end{equation}
where
\begin{equation}
M_{2D}=%
\begin{pmatrix}
I_{D} & t\,m^{-1}\\
0 & I_{D}%
\end{pmatrix}
\label{eq:Stability}%
\end{equation}
is the stability matrix and $M_{2D^{2}}=M_{2D}\otimes I_{D}$. The kinetic flow
with Jacobian~(\ref{eq:HJT}) is symplectic, i.e., $\Phi_{\text{T},t}^{\prime
}(z_{0})^{T}\cdot\omega(z_{t})\cdot\Phi_{\text{T},t}^{\prime}(z_{0}%
)=\omega(z_{0}),$ because
\begin{equation}
M_{2D}^{T}\cdot J_{2D}\cdot M_{2D}=%
\begin{pmatrix}
I_{D} & 0\\
t\,m^{-1} & I_{D}%
\end{pmatrix}%
\begin{pmatrix}
0 & I_{D}\\
-I_{D} & -t\,m^{-1}%
\end{pmatrix}
=J_{2D}%
\end{equation}
and, due to the relation between matrix and tensor multiplications,
\begin{align}
M_{2D^{2}}^{T}  &  \cdot J_{2D^{2}}\cdot M_{2D^{2}}\nonumber\\
&  =\big(M_{2D}^{T}\otimes I_{D}\big)\cdot\big(J_{2D}\otimes I_{D}%
\big)\cdot\big(M_{2D}\otimes I_{D}\big)\nonumber\\
&  =\big(M_{2D}^{T}\cdot J_{2D}\cdot M_{2D}\big)\otimes I_{D}=J_{2D^{2}}.
\label{eq:Stability_is_symplectic_2}%
\end{align}


\subsection{Symplecticity of the potential flow}

Inserting the coefficients $V_{0}$, $V_{1}$, and $V_{2}$ of the fully
variational GWD [Eq.~\eqref{eq:V0_V1_V2_VGA}], the local harmonic GWD
[Eq.~\eqref{eq:V0_V1_V2_TGA}], or the local cubic variational GWD
[Eq.~\eqref{eq:V0_V1_V2_LCA}] into Eqs.~(\ref{eq:Hf_V_qEOM_LCA}%
)-(\ref{eq:Hf_V_PEOM_LCA}) gives the potential flow $\Phi_{\text{V},t}(z_{0})$
for the corresponding methods. The Jacobian of this potential flow is the
$(2D+4D^{2})$-dimensional matrix
\begin{equation}
\Phi_{\text{V},t}^{\prime}(z_{0})=I_{2D+4D^{2}}-t%
\begin{pmatrix}
0 & 0 & 0 & 0 & 0 & 0\\
a & 0 & b^{(1)} & 0 & b^{(2)} & 0\\
0 & 0 & 0 & 0 & 0 & 0\\
c^{(1)} & 0 & d^{(11)} & 0 & d^{(12)} & 0\\
0 & 0 & 0 & 0 & 0 & 0\\
c^{(2)} & 0 & d^{(21)} & 0 & d^{(22)} & 0
\end{pmatrix}
, \label{eq:HJV}%
\end{equation}
where
\begin{align}
a_{j,k}  &  =\partial V_{1}(q_{0},Q_{0})_{j}/\partial q_{0,k},\label{eq:a}\\
b_{j,D(k-1)+l}^{(r)}  &  =\partial V_{1}(q_{0},Q_{0})_{j}/\partial
Q_{0,kl}^{(r)},\label{eq:b}\\
c_{D(j-1)+k,l}^{(r)}  &  =\partial\big [V_{2}(q_{0},Q_{0})_{jm}\,Q_{0,mk}%
^{(r)}\big]/\partial q_{0,l},\label{eq:c}\\
d_{D(j-1)+k,\,D(l-1)+m}^{(rs)}  &  =\partial\big [V_{2}(q_{0},Q_{0}%
)_{jn}\,Q_{0,nk}^{(r)}\big]/\partial Q_{0,lm}^{(s)}\nonumber\\
&  =\big [\partial V_{2}(q_{0},Q_{0})_{jn}/\partial Q_{0,lm}^{(s)}%
\big]\,Q_{0,nk}^{(r)}\nonumber\\
&  \quad+V_{2}(q_{0},Q_{0})_{jn}\,\delta_{nl}\,\delta_{km}\,\delta_{rs}
\label{eq:d}%
\end{align}
are, respectively, components of a $D\times D$ matrix $a$, a $D\times D^{2}$
matrix $b^{(r)}$, a $D^{2}\times D$ matrix $c^{(r)}$, and a $D^{2}\times
D^{2}$ matrix $d^{(rs)}$ for all $r,s\in\{1,2\}$. Satisfaction of the relation
$\Phi_{\text{V},t}^{\prime}(z_{0})^{T}\cdot\omega(z_{t})\cdot\Phi_{\text{V}%
,t}^{\prime}(z_{0})=\omega(z_{0})$ requires that
\begin{equation}
a^{T}=a,\qquad\big(b^{(r)}\big)^{T}=(\hbar/2)\,c^{(r)},\qquad\big(d^{(rs)}%
\big)^{T}=d^{(sr)}. \label{eq:conditions}%
\end{equation}

Since the coefficients~\eqref{eq:V0_V1_V2_VGA} of the fully variational GWD
depend on $\langle\hat{V}\rangle$, $\langle\hat{V}^{\prime}\rangle$, and
$\langle\hat{V}^{\prime}\rangle$, which, in turn, depend only on $q_{0}$ and
$Q_{0}$, in Eqs.~(\ref{eq:a})-(\ref{eq:d}) we only need expressions for the
derivatives of $\langle\hat{V}^{(n)}\rangle$ with respect to $q_{0}$ and
$Q_{0}$. The former is
\begin{equation}
\partial\langle\hat{V}^{(n)}\rangle/\partial q_{0}=\langle\hat{V}%
^{(n+1)}\rangle,
\end{equation}
while the latter can be found in Eq.~(B30) in
Ref.~\onlinecite{Fereidani_Vanicek:2023}. Using these two relations, we
conclude that
\begin{align}
a_{j,k} &  =\langle\hat{V}^{\prime\prime}\rangle_{jk},\label{eq:a_VGA}\\
b_{j,D(k-1)+l}^{(r)} &  =(\hbar/2)\,\langle\hat{V}^{\prime\prime\prime}%
\rangle_{jkm}\,Q_{0,ml}^{(r)},\label{eq:b_VGA}\\
c_{D(j-1)+k,l}^{(r)} &  =\langle\hat{V}^{\prime\prime\prime}\rangle
_{jml}\,Q_{0,mk}^{(r)},\label{eq:c_VGA}\\
d_{D(j-1)+k,\,D(l-1)+m}^{(rs)} &  =(\hbar/2)\,\langle\hat{V}^{(4)}%
\rangle_{jnlp}\,Q_{0,nk}^{(r)}\,Q_{0,pm}^{(s)}\nonumber\\
&  \quad+\langle\hat{V}^{\prime\prime}\rangle_{jn}\,\delta_{nl}\,\delta
_{km}\,\delta_{rs}.\label{eq:d_VGA}%
\end{align}
Because $\langle\hat{V}^{\prime\prime}\rangle$, $\langle\hat{V}^{\prime
\prime\prime}\rangle$, and $\langle\hat{V}^{(4)}\rangle$ are totally
symmetric, Eqs.~(\ref{eq:a_VGA})-(\ref{eq:d_VGA}) for $a$, $b^{(r)}$,
$c^{(r)}$, and $d^{(rs)}$ imply that conditions~(\ref{eq:conditions}) hold for
the Jacobian~(\ref{eq:HJV}), and thus the potential flow of the variational
GWD is symplectic.

Finally, inserting coefficients~\eqref{eq:V0_V1_V2_TGA} of the local harmonic
GWD into Eqs.~(\ref{eq:a})-(\ref{eq:d}) yields
\begin{align}
a_{j,k}  &  =V^{\prime\prime}(q_{0})_{jk},\label{eq:a_TGA}\\
b_{j,D(k-1)+l}^{(r)}  &  =0,\label{eq:b_TGA}\\
c_{D(j-1)+k,l}^{(r)}  &  =V^{\prime\prime\prime}(q_{0})_{jml}\,Q_{0,mk}%
^{(r)},\label{eq:c_TGA}\\
d_{D(j-1)+k,\,D(l-1)+m}^{(rs)}  &  =V^{\prime\prime}(q_{0})_{jn}\,\delta
_{nl}\,\delta_{km}\,\delta_{rs}. \label{eq:d_TGA}%
\end{align}
Equations~(\ref{eq:b_TGA}) and~(\ref{eq:c_TGA}) imply that the second
condition in~(\ref{eq:conditions}) is not fulfilled, and therefore the local
harmonic GWD is not symplectic.

Inserting coefficients~\eqref{eq:V0_V1_V2_LCA} of the local cubic variational
GWD into Eqs.~(\ref{eq:a})-(\ref{eq:d}) yields
\begin{align}
a_{j,k}  &  =V^{\prime\prime}(q_{0})_{jk}\nonumber\\
&  \quad+V^{(4)}(q_{0})_{jklm}\,\Sigma_{0,lm}/2,\label{eq:a_LCA}\\
b_{j,D(k-1)+l}^{(r)}  &  =(\hbar/2)\,V^{\prime\prime\prime}(q_{0}%
)_{jkm}\,Q_{0,ml}^{(r)},\label{eq:b_LCA}\\
c_{D(j-1)+k,l}^{(r)}  &  =V^{\prime\prime\prime}(q_{0})_{jml}\,Q_{0,mk}%
^{(r)},\label{eq:c_LCA}\\
d_{D(j-1)+k,\,D(l-1)+m}^{(rs)}  &  =V^{\prime\prime}(q_{0})_{jn}\,\delta
_{nl}\,\delta_{km}\,\delta_{rs}, \label{eq:d_LCA}%
\end{align}
where we used the position covariance in Hagedorn's
parametrization,~\cite{Vanicek_Begusic:2021} $\Sigma_{t}:=(\hbar
/2)\,Q_{t}\cdot Q_{t}^{\dagger}$, to derive~\eqref{eq:b_LCA}. The
conditions~(\ref{eq:conditions}) are satisfied for~(\ref{eq:HJV}) because
$V^{\prime\prime}$, $V^{\prime\prime\prime}$, and the fourth derivative
$V^{(4)}$ of the potential, in Eqs.(\ref{eq:a_LCA})-(\ref{eq:d_LCA}), are
totally symmetric. Hence, the potential flow of the local cubic variational
GWD is symplectic.

\subsection{Symplecticity of the composed geometric integrators}

If both kinetic and potential flows are symplectic, then any composition of
them is symplectic. This proves the conservation of the symplectic
structure~(\ref{eq:sympl_QPS}) by the geometric integrators developed for the
fully variational and local cubic variational GWD, since the composed flow
$\Phi_{t}(z_{0})$ of a geometric integrator consists of many steps, each of
which consists of several potential and kinetic substeps. In
Ref.~\onlinecite{Fereidani_Vanicek:2023} and in
Sec.~\ref{subsec:2D_coupled_Morse}, we verified symplecticity of these
integrators numerically by measuring the accuracy with which
Eq.~(\ref{eq:JBJ=B}) is satisfied if $\Phi_{t}(z_{0})$ denotes the composed
flow consisting many steps, each of which, in turn, is composed by several
kinetic and potential substeps. For that, Jacobian $\Phi_{t}^{\prime}(z_{0})$
of the composed flow appearing in Eq.~(\ref{eq:JBJ=B}) was obtained by matrix
multiplication of the Jacobians of all kinetic and potential steps.


\bibliographystyle{aipnum4-2}
\bibliography{biblio53,addition_VGA_LCA,duplicates_VGA_LCA}

\end{document}